\documentclass[onefignum,onetabnum]{siamart251216}

\usepackage{amsfonts}
\usepackage{mathrsfs}
\usepackage{amsmath}
\usepackage{amssymb}
\usepackage{mathtools}
\usepackage{algorithm}
\usepackage{algpseudocode}

\headers{Learned Adaptive Multiresolution Diffusion Imaging}
{C. Tantardini et al.}

\title{Learned Adaptive Multiresolution Diffusion Imaging
}

\author{
Christian Tantardini
\thanks{Center for Integrative Petroleum Research,
King Fahd University of Petroleum and Minerals,
Dhahran 31261, Kingdom of Saudi Arabia
(\email{christiantantardini@ymail.com}).
}
\and
Stig Rune Jensen
\thanks{
Hylleraas Centre for Quantum Molecular Sciences,
Department of Chemistry,
UiT The Arctic University of Norway,
N-9037 Tromsø, Norway
(\email{stig.r.jensen@uit.no}).
}
\and
Roberto Di Remigio Eik{\AA}s
\thanks{Algorithmiq S.r.l., c/o STLex, Via della Chiusa 15, 20123, Milano, Italia
(\email{roberto@algorithmiq.tech})
}
\and
Joakim Henrik Beck
\thanks{
Center for Integrative Petroleum Research,
King Fahd University of Petroleum and Minerals,
Dhahran 31261, Kingdom of Saudi Arabia.
(\email{joakim.beck@kfupm.edu.sa}).
}
}

\usepackage{hyperref}
\usepackage{cleveref}

\newcommand{\Energy}{\mathcal{E}}

\begin{document}

\maketitle

\begin{abstract}
Adaptive multiresolution methods reduce representation cost by concentrating fine-scale degrees of freedom where needed, but their tree updates are usually governed by fixed local criteria. We introduce Learned Adaptive Multiresolution Diffusion Imaging (Learned AMDI), which preserves the AMDI fixed-tree propagator and hierarchy constraints while replacing the post-propagation selector with a shared local policy trained by proximal policy optimization. Regression tests reproduce deterministic AMDI trajectories to machine precision when identical trees are used. In the Haar implementation studied here, the deterministic one-step selector accepts no refinements in 54 decisions. Across nine held-out cases, Learned AMDI executes 393 refinements and reduces the mean terminal reference discrepancy from $0.17496$ to $0.13657$, while occupancy rises from $0.13737$ to $0.26660$. Step-resolved diagnostics reveal occasional small adaptation-energy increases; fixed-tree energy stability therefore does not guarantee monotonicity of the learned outer iteration. At comparable occupancy, a validation-tuned observed-detail threshold reaches a discrepancy of $0.13792$ with slightly better RMSE and SSIM, placing both methods on essentially the same accuracy--occupancy tradeoff. A decision-1-only control reaches $0.13742$, indicating that most of the improvement on this static benchmark arises from the initial allocation. The shared actor transfers without retraining to $64\times64$ and $128\times128$ images, improving reference discrepancy, RMSE, and SSIM relative to deterministic AMDI, while the frozen threshold rule remains competitive. Learned AMDI thus provides a hierarchy-constrained, resolution-transferable mechanism for adaptive allocation and clarifies the contribution of sequential decisions.
\end{abstract}

\begin{keywords}
adaptive multiresolution methods, reinforcement learning, diffusion imaging,
wavelets, adaptive discretization, proximal policy optimization
\end{keywords}

\begin{MSCcodes}
65N12, 65N30, 65T60, 49J40
\end{MSCcodes}

%%%%%%%%%%%%%%%%%%%%%%%%%%%%%%%%%%%%%%%%%%%%%%%%%%%%%%%%%%%%%%%%%%%%%%%%%%
\section{Introduction}
\label{sec:introduction}
%%%%%%%%%%%%%%%%%%%%%%%%%%%%%%%%%%%%%%%%%%%%%%%%%%%%%%%%%%%%%%%%%%%%%%%%%%

Diffusion-based methods remain attractive for image denoising because their
smoothing mechanisms follow from explicit models
\cite{PeronaMalik1990,CatteEtAl1992,Weickert1999,BlackEtAl1998}.
Total-variation and nonlocal formulations provide complementary
edge-preserving viewpoints
\cite{RudinOsherFatemi1992,Chambolle2004,BuadesCollMorel2005}.
Their numerical representations, however, need not have uniform resolution.
Smooth regions may require few degrees of freedom, whereas edges and localized
textures require finer spatial detail. Wavelet and multiresolution frameworks
exploit this disparity by retaining fine-scale information only where it is
needed
\cite{Mallat1989,Daubechies1992,Harten1995,Harten1996,DeVore1998,CohenDahmenDeVore2001,CohenEtAl2003}.
Data-dependent wavelet shrinkage provides a closely related precedent for
spatially selective denoising
\cite{DonohoJohnstone1994,Donoho1995}; adaptive representations have also been
used directly in diffusion-based imaging \cite{PreusserRumpf2000}.

The previous work on deterministic adaptive multiresolution diffusion imaging (AMDI) \cite{TantardiniAMDI2026} developed a variational method in which the active multiresolution tree, its coefficients, and the associated interaction graph together form the adaptive state. Each outer iteration combines a fixed-tree coefficient update with a deterministic choice among admissible tree modifications. Because the interaction operator depends on the active hierarchy, changing the tree also changes the graph used in subsequent updates. In the present work, we preserve the AMDI coefficient update, numerical parameters, and hierarchy constraints, and replace only the rule used to select the next tree.

This distinction is important because a tree modification can have a delayed benefit. Refinement may increase the immediate energy and representation cost while making additional coefficients available to later propagation steps. A selector based only on the current penalized energy can therefore reject a refinement that would improve the terminal state. We address this limitation using a learned policy trained from complete trajectories, while explicitly testing whether repeated decisions provide a measurable advantage over a single early adaptation.

Classical adaptive mesh refinement established hierarchical local resolution as
a general strategy for balancing discretization error and computational work
\cite{BergerOliger1984,BergerColella1989}. More recently, reinforcement
learning has been used to select adaptive meshes during numerical simulations
\cite{YangEtAl2023,FoucartEtAl2023,FreymuthEtAl2023,DzanicEtAl2024}.
These studies provide precedents for policies that act on local information,
optimize trajectory-level objectives, and transfer across problem sizes.
Within scientific machine learning, physics-informed models incorporate the
governing equations into training \cite{RaissiEtAl2019,KarniadakisEtAl2021},
whereas neural operators and multiscale operator architectures learn solution
maps between function spaces
\cite{LuEtAl2021DeepONet,LiEtAl2021FNO,KovachkiEtAl2023,GuptaEtAl2021,PandeyEtAl2026LiNO}.
Graph-based simulators provide a further route to discretization-aware learned
dynamics \cite{PfaffEtAl2021,BrandstetterEtAl2022}. Our objective is narrower:
the learned policy proposes admissible modifications of the active AMDI tree,
while the image field continues to evolve through the prescribed numerical
operator.

Let $\left(u^{n},\mathcal{M}^{n}\right)$ denote the accepted field and multiresolution representation after adaptation at time $t_n$. The next outer iteration first propagates the field while holding the representation fixed:
\begin{equation}
\widetilde{u}^{\,n+1}
=
\mathcal{D}_{\Delta t}
\left(u^{n};\mathcal{M}^{n}\right).
\label{eq:intro_propagation}
\end{equation}
The policy then acts on a state $\mathcal{S}^{n}$ constructed from the propagated field, the observed image, and the current hierarchy:
\begin{equation}
a^{n}
\sim
\pi_{\theta}\left(\,\cdot\,\middle|\mathcal{S}^{n}\right).
\label{eq:intro_policy}
\end{equation}
The local actions request refinement, retention, or coarsening. A deterministic multiresolution update resolves these requests and enforces hierarchical admissibility before the next coefficient update.

Deterministic and Learned AMDI therefore use the same fixed-tree propagation and the same propagate--then--adapt ordering. Their trajectories differ only through the selected trees, which produce the active coefficients and subsequent interaction graphs. The energy inequality established for the variational selector in the companion work does not automatically extend to a learned tree update. We consequently record the energy immediately before and after every adaptation, rather than assuming monotonic decay of the learned outer iteration.

The shared local policy is trained using proximal policy optimization. Its trajectory reward combines discrepancy from a uniformly resolved AMDI calculation, representation occupancy, and hierarchy switching. The uniformly resolved trajectory is a numerical reference rather than the clean image; reference discrepancy is therefore reported separately from RMSE and SSIM against the clean target. Similarly, occupancy measures the fraction of active degrees of freedom and should not be interpreted directly as execution time or memory consumption.

The numerical study is designed to identify what the selectors actually do. On the same nine holdout cases and noise realizations, we record proposed and executed actions, active degrees of freedom, deterministic candidate energies, and pre- and post-adaptation energy. We compare Learned AMDI not only with deterministic AMDI, but also with retention of the initial tree and two independently validation-tuned observed-detail controls. A further temporal control permits the actor to act only at the first decision.

The numerical results identify the role of the learned selector in the implemented AMDI configuration.
Zero initialization creates a structural barrier for the deterministic one-step refinement candidate, and the deterministic selector consequently accepts no refinements in the 54 holdout decisions.
Learned AMDI executes 393 refinements and reduces the mean terminal reference error from \(0.17496\) to \(0.13657\), a reduction of \(21.9\%\), while increasing the active representation.
At nearby but not identical occupancy, the validation-tuned threshold rule reaches essentially the same accuracy--occupancy tradeoff and gives slightly better RMSE and SSIM.
The principal demonstrated benefit is therefore objective-driven allocation under the AMDI hierarchy constraints.
The temporal control shows that this allocation is formed almost entirely at the first decision on the static benchmark, while the shared local actor transfers without retraining to \(64\times64\) and \(128\times128\) images.

The contribution of this work is a finite-horizon formulation of AMDI tree selection, a shared local policy that preserves hierarchical admissibility, and a step-resolved numerical analysis connecting proposed and executed actions to representation occupancy, reference accuracy, and energy change.
Together, these elements provide a transparent learned mechanism for selecting an adaptive representation while retaining the established AMDI propagation procedure.

Section~\ref{sec:amdi_background} recalls the AMDI construction.
Section~\ref{sec:sequential_formulation} defines the decision process,
and Section~\ref{sec:learning_method} describes the policy and its
training. Section~\ref{sec:numerical_experiments} presents the
numerical protocol, controls, diagnostics, and transfer experiments.
Section~\ref{sec:conclusions} summarizes the conclusions and
limitations.

%%%%%%%%%%%%%%%%%%%%%%%%%%%%%%%%%%%%%%%%%%%%%%%%%%%%%%%%%%%%%%%%%%%%%%%%%%
\section{Deterministic Adaptive Multiresolution Diffusion Imaging}
\label{sec:amdi_background}
%%%%%%%%%%%%%%%%%%%%%%%%%%%%%%%%%%%%%%%%%%%%%%%%%%%%%%%%%%%%%%%%%%%%%%%%%%

We first summarize the components of the deterministic Adaptive
Multiresolution Diffusion Imaging (AMDI) formulation that are required for
the learned extension developed below. The purpose of this section is not to
repeat the complete derivation of AMDI, which is given in the companion
paper~\cite{TantardiniAMDI2026}, but to establish the notation and identify
precisely which part of the original algorithm is replaced by learning.

%%%%%%%%%%%%%%%%%%%%%%%%%%%%%%%%%%%%%%%%%%%%%%%%%%%%%%%%%%%%%%%%%%%%%%%%%%
\subsection{Diffusion evolution}
\label{sec:diffusion_evolution}
%%%%%%%%%%%%%%%%%%%%%%%%%%%%%%%%%%%%%%%%%%%%%%%%%%%%%%%%%%%%%%%%%%%%%%%%%%

Let $\Omega\subset\mathbb{R}^{d}$ be the spatial domain. In the
denoising experiments $d=2$, and $u^n$ denotes the image field after
the $n$th accepted adaptive update. For an active representation
$\mathcal{M}^{n}$, write $V(\mathcal{M}^{n})$ for its approximation
space. The next field update is performed while this representation
is held fixed:
\begin{equation}
    \widetilde{u}^{\,n+1}
    =
    \mathcal{D}_{\Delta t}
    \left(u^n;\mathcal{M}^{n}\right).
    \label{eq:discrete_diffusion_operator}
\end{equation}

In the variational AMDI formulation, the corresponding ideal
fixed-tree step minimizes an incremental term together with the
AMDI energy:
\begin{equation}
    \widetilde{u}^{\,n+1}_{\mathrm{ideal}}
    \in
    \operatorname*{arg\,min}_{v\in V(\mathcal{M}^{n})}
    \left\{
        \frac{1}{2\Delta t}\|v-u^n\|_2^2
        +
        \Energy(\mathcal{M}^{n};v)
    \right\}.
    \label{eq:ideal_fixed_tree_step}
\end{equation}

Writing \(c=(c_i)\) for the active Haar coefficients and \(u_{\mathcal{M},c}\) for their reconstruction, the implemented AMDI energy is
\begin{equation}
\begin{aligned}
\Energy(\mathcal{M};c)
={}&
\frac{1}{2}
\left\|
u_{\mathcal{M},c}-u_{\mathrm{obs}}
\right\|_{2}^{2}
+
\frac{\alpha}{4}
\sum_{i,j}
w_{ij}(\mathcal{M},c)
(c_i-c_j)^2
\\
&+
\beta\sum_i \nu_i |c_i|
+
\tau\sum_i
\left(1+\gamma_{\mathrm{lev}}\ell_i\right),
\qquad
\nu_i=1+0.05\ell_i .
\end{aligned}
\label{eq:implemented_amdi_energy}
\end{equation}
Here \(\ell_i\) is the level of coefficient \(i\), and the symmetric weights \(w_{ij}\geq0\) define the state-dependent interaction graph. The factor \(\alpha/4\) corresponds to summation over ordered index pairs. The four terms represent data fidelity, graph interaction, coefficient sparsity, and tree complexity, respectively. During a fixed-tree minimization, the final term is constant.

The reported Haar implementation uses a frozen-interaction
approximation: interaction weights are assembled from the accepted
state, the resulting coefficient problem is solved approximately,
and a backtracking check tests the original AMDI energy before the
coefficient update is accepted. This safeguard concerns the
fixed-tree update. Tree selection takes place afterward. Thus,
$\mathcal{D}_{\Delta t}$ denotes the same numerical procedure in
both methods, while its assembled interactions depend on the
representation and coefficients supplied to it.

%%%%%%%%%%%%%%%%%%%%%%%%%%%%%%%%%%%%%%%%%%%%%%%%%%%%%%%%%%%%%%%%%%%%%%%%%%
\subsection{Multiresolution representation}
\label{sec:multiresolution_representation}
%%%%%%%%%%%%%%%%%%%%%%%%%%%%%%%%%%%%%%%%%%%%%%%%%%%%%%%%%%%%%%%%%%%%%%%%%%

At each diffusion step, AMDI represents the field on an adaptive
multiresolution structure rather than on a permanently active finest-level
discretization.
We denote the complete active representation at step $n$ by
$\mathcal{M}^{n}$.
The object $\mathcal{M}^{n}$ contains the active multiresolution degrees of
freedom required to represent $u^{n}$ according to the adaptive hierarchy of
deterministic AMDI~\cite{TantardiniAMDI2026}.

This representation inherits the hierarchical character of standard
multiresolution methods, in which spatial information is organized across
resolution levels and localized fine-scale contributions can be introduced or
removed according to their relevance
\cite{Mallat1989,Harten1996,CohenEtAl2003}. We intentionally retain the
compact notation $\mathcal{M}^{n}$ here because the detailed construction of
the deterministic hierarchy is given in the companion paper and is not
modified in the present work.

Let $N_{\mathrm{act}}^{n}$ denote the number of active multiresolution
degrees of freedom contained in $\mathcal{M}^{n}$, and let
$N_{\mathrm{full}}$ denote the corresponding number for the uniformly
resolved finest representation.
We define the
instantaneous representation occupancy as
\begin{equation}
    \rho_{\mathrm{act}}^{n}
    =
    \frac{
        N_{\mathrm{act}}^{n}
    }{
        N_{\mathrm{full}}
    },
    \qquad
    0 < \rho_{\mathrm{act}}^{n}\leq 1.
    \label{eq:representation_occupancy}
\end{equation}
Thus, $\rho_{\mathrm{act}}^{n}=1$ corresponds to the complete uniformly
resolved representation, whereas $\rho_{\mathrm{act}}^{n}\ll1$ indicates
that only a small fraction of the available multiresolution degrees of
freedom is active.

We distinguish representation occupancy from implementation-dependent
computational cost.  In the present work,
$\rho_{\mathrm{act}}^{n}$ and its terminal value
$C_{\mathrm{rel}}$ used in the numerical summaries are therefore measures of
representation complexity rather than direct surrogates for wall-clock time
or memory consumption.

%%%%%%%%%%%%%%%%%%%%%%%%%%%%%%%%%%%%%%%%%%%%%%%%%%%%%%%%%%%%%%%%%%%%%%%%%%
\subsection{Deterministic adaptive evolution}
\label{sec:deterministic_amdi}
%%%%%%%%%%%%%%%%%%%%%%%%%%%%%%%%%%%%%%%%%%%%%%%%%%%%%%%%%%%%%%%%%%%%%%%%%%

The deterministic AMDI algorithm used as the baseline in the present work is
the adaptive diffusion procedure introduced in
Ref.~\cite{TantardiniAMDI2026}.  Its outer iteration follows a
propagate--then--adapt sequence.  Starting from the accepted field
$u^{n}$ represented on the active hierarchy $\mathcal{M}^{n}$, the field is
first advanced over one diffusion increment while the hierarchy is held
fixed:
\begin{equation}
    \widetilde{u}^{\,n+1}
    =
    \mathcal{D}_{\Delta t}
    \left(
        u^{n};
        \mathcal{M}^{n}
    \right).
    \label{eq:amdi_deterministic_propagation}
\end{equation}
Here $\widetilde{u}^{\,n+1}$ denotes the propagated field before the adaptive
representation is updated.

The deterministic AMDI adaptation is then applied to the propagated field.
For compactness, we introduce in the present work the notation
$\mathcal{U}_{\mathrm{det}}$ for this complete deterministic adaptive update:
\begin{equation}
    \left(
        u^{n+1},
        \mathcal{M}^{n+1}
    \right)
    =
    \mathcal{U}_{\mathrm{det}}
    \left(
        \widetilde{u}^{\,n+1},
        \mathcal{M}^{n}
    \right).
    \label{eq:deterministic_amdi_update}
\end{equation}
The symbol $\mathcal{U}_{\mathrm{det}}$ is used here only as an abstract
representation of the deterministic adaptation procedure of
Ref.~\cite{TantardiniAMDI2026}; it is not an additional operator introduced
into the underlying AMDI dynamics.

For the denoising implementation used below,
$\mathcal{U}_{\mathrm{det}}$ constructs at most three candidate
trees: the current tree, a tree formed by refining the
highest-scoring admissible leaves, and a tree formed by coarsening
the lowest-scoring prunable groups. The scores use Haar details of
the observed image; the reported refinement and coarsening
fractions are specified in Section~\ref{sec:numerical_experiments}.
The propagated coefficients are transferred to each candidate
representation, and newly opened detail coefficients are
initialized to zero in the reported implementation. The selected
candidate minimizes the AMDI energy after transfer plus the
prescribed tree-change penalty
\cite{TantardiniAMDI2026}.

Because the unchanged tree is among the candidates, exact
minimization of the tree-selection criterion gives a selected
candidate whose energy plus tree-change penalty is no larger
than that of the post-propagation state on the unchanged tree.

For the reported Haar implementation, zero initialization creates a
structural refinement barrier.
Transferring the propagated coefficients to a strictly refined tree while setting every newly opened detail coefficient to zero leaves the reconstructed field, and therefore the fidelity term in
Eq.~\eqref{eq:implemented_amdi_energy}, unchanged.
It also leaves the nonsmoothed sparsity term unchanged.
For the refinement-stable interaction graph used here, refinement preserves all pre-existing interactions and can add only nonnegative quadratic contributions.
At the same time, the positive tree-complexity term strictly increases.
The nonnegative tree-change penalty can only increase this difference.
Consequently, every strict refinement candidate has a larger implemented selection objective than the unchanged-tree candidate and cannot be selected by \(\mathcal{U}_{\mathrm{det}}\).
This property is specific to the implemented transfer-and-score construction; it need not hold if newly opened details are initialized from the observed data or if candidate coefficients are refitted before the trees are compared. 
Coarsening and retention remain data dependent.

The representation selected at step $n$ is used in the next
fixed-tree update. Both AMDI variants use this
propagate--then--adapt ordering.

%%%%%%%%%%%%%%%%%%%%%%%%%%%%%%%%%%%%%%%%%%%%%%%%%%%%%%%%%%%%%%%%%%%%%%%%%%
\subsection{From deterministic to learned adaptive evolution}
\label{sec:prescribed_to_learned}
%%%%%%%%%%%%%%%%%%%%%%%%%%%%%%%%%%%%%%%%%%%%%%%%%%%%%%%%%%%%%%%%%%%%%%%%%%

Deterministic AMDI and Learned AMDI differ only in the mechanism used to obtain the
adaptive representation after propagation.  In both methods, the accepted
state
$\left(u^{n},\mathcal{M}^{n}\right)$ is first advanced using the same
diffusion propagator,
\begin{equation}
    \widetilde{u}^{\,n+1}
    =
    \mathcal{D}_{\Delta t}
    \left(
        u^{n};
        \mathcal{M}^{n}
    \right).
    \label{eq:common_predecision_propagation}
\end{equation}

For deterministic AMDI, the propagated field is processed by the
deterministic adaptive procedure summarized by
$\mathcal{U}_{\mathrm{det}}$,
\begin{equation}
    \left(
        u^{n+1},
        \mathcal{M}^{n+1}
    \right)
    =
    \mathcal{U}_{\mathrm{det}}
    \left(
        \widetilde{u}^{\,n+1},
        \mathcal{M}^{n}
    \right).
    \label{eq:deterministic_postprop_adaptation}
\end{equation}

In Learned AMDI, the adaptive selection is instead obtained by a
parameterized policy.  The policy acts on information constructed from the
propagated field $\widetilde{u}^{\,n+1}$ and the current hierarchy
$\mathcal{M}^{n}$,
\begin{equation}
    a^{n}
    \sim
    \pi_{\theta}
    \left(
        \,\cdot\,
        \middle|
        \mathcal{S}^{n}
    \right),
    \label{eq:learned_policy}
\end{equation}
and the selected actions are applied through the constrained
multiresolution-update operator,
\begin{equation}
    \left(
        u^{n+1},
        \mathcal{M}^{n+1}
    \right)
    =
    \mathcal{U}_{\mathrm{MR}}
    \left(
        \widetilde{u}^{\,n+1},
        \mathcal{M}^{n},
        a^{n}
    \right).
    \label{eq:learned_representation_update}
\end{equation}

For identical accepted inputs, the two methods call the same
fixed-tree update and obtain the same propagated field. They then
choose the next representation differently:
$\mathcal{U}_{\mathrm{det}}$ evaluates its prescribed candidate
trees, whereas $\mathcal{U}_{\mathrm{MR}}$ resolves the policy's
local proposals subject to hierarchy constraints. Once their
representations differ, subsequent interaction graphs and
propagated fields can differ as well.

The learned update preserves the admissibility of the resulting
hierarchy. It is not, however, the energy-minimizing tree update
used in the companion AMDI energy argument. That argument does
not automatically give energy decay for a learned outer
iteration.

%%%%%%%%%%%%%%%%%%%%%%%%%%%%%%%%%%%%%%%%%%%%%%%%%%%%%%%%%%%%%%%%%%%%%%%%%%
\section{Adaptive Multiresolution Diffusion as a Sequential Decision Problem}
\label{sec:sequential_formulation}
%%%%%%%%%%%%%%%%%%%%%%%%%%%%%%%%%%%%%%%%%%%%%%%%%%%%%%%%%%%%%%%%%%%%%%%%%%

The deterministic AMDI formulation updates the multiresolution
representation through the adaptive procedure summarized in
Section~\ref{sec:deterministic_amdi}.  Learned AMDI replaces the
deterministic selection mechanism within this post-propagation adaptive stage
by a policy optimized from complete diffusion trajectories.  Because a
decision made at step $n$ modifies the representation available to all
subsequent diffusion updates, the adaptive problem is inherently sequential.
We therefore formulate Learned AMDI as a finite-horizon Markov decision
process (MDP)
\cite{KaelblingEtAl1996,Puterman1994}, following the general viewpoint used in
reinforcement-learning approaches to adaptive numerical discretization
\cite{YangEtAl2023,FoucartEtAl2023}.

We denote the finite-horizon decision process by
\begin{equation}
    \mathfrak{M}
    =
    \left(
        \mathscr{S},
        \mathscr{A},
        \mathcal{P},
        R,
        \mu_{0},
        \gamma,
        N_{t}
    \right),
    \label{eq:mdp_definition}
\end{equation}
where $\mathscr{S}$ is the state space,
$\mathscr{A}$ is the adaptive action space,
$\mathcal{P}$ is the transition kernel,
$R$ is the stage-reward function,
$\mu_{0}$ is the initial-state distribution,
$\gamma\in(0,1]$ is the discount factor, and
$N_{t}$ is the finite number of diffusion steps in one episode.
The individual components are defined below so that the learned decision
process remains compatible with the deterministic AMDI dynamics introduced
in Section~\ref{sec:amdi_background}.

%%%%%%%%%%%%%%%%%%%%%%%%%%%%%%%%%%%%%%%%%%%%%%%%%%%%%%%%%%%%%%%%%%%%%%%%%%
\subsection{State representation}
\label{sec:rl_state}
%%%%%%%%%%%%%%%%%%%%%%%%%%%%%%%%%%%%%%%%%%%%%%%%%%%%%%%%%%%%%%%%%%%%%%%%%%

It is useful to distinguish the accepted AMDI state from the state at which
the adaptive decision is made.  After adaptation at time $t_n$, let
\begin{equation}
    \mathcal{X}^{n}
    =
    \left(
        u^{n},
        \mathcal{M}^{n},
        t_n
    \right)
    \label{eq:accepted_amdi_state}
\end{equation}
denote the accepted physical--multiresolution state.

Consistently with the deterministic AMDI outer iteration, the next decision
is not made from $\mathcal{X}^{n}$ directly.  The field is first propagated on the
current representation:
\begin{equation}
    \widetilde{u}^{\,n+1}
    =
    \mathcal{D}_{\Delta t}
    \left(
        u^{n};
        \mathcal{M}^{n}
    \right).
    \label{eq:decision_state_propagation}
\end{equation}

The features supplied to the policy are then constructed from this
post-propagation, pre-adaptation field,
\begin{equation}
    \boldsymbol{\phi}^{n}
    =
    \mathcal{F}
    \left(
        u_{\mathrm{obs}},
        \widetilde{u}^{\,n+1},
        \mathcal{M}^{n},
        t_{n+1}
    \right).
    \label{eq:state_feature_map}
\end{equation}

We therefore define the MDP decision state at outer iteration $n$ as
\begin{equation}
    \mathcal{S}^{n}
    =
    \left(
        u_{\mathrm{obs}},
        \widetilde{u}^{\,n+1},
        \mathcal{M}^{n},
        \boldsymbol{\phi}^{n},
        t_{n+1}
    \right)
    \in
    \mathscr{S}.
    \label{eq:rl_state_definition}
\end{equation}
Here $u_{\mathrm{obs}}$ is fixed throughout one episode and enters the
environment state because part of the local feature construction depends on
the observed image.  The actor does not receive the complete observed image
directly; it receives only the local feature vectors defined in
Section~\ref{sec:policy_features}.  The index $n$ labels the adaptive
decision stage, while the propagated field contained in
$\mathcal{S}^{n}$ has already reached $t_{n+1}$.

For fixed numerical parameters and a specified observed image,
the full decision state in
Eq.~\eqref{eq:rl_state_definition} contains the information needed
to get the next state after an action is applied. The same
fixed-tree update routine is then called on the newly accepted
representation. In this sense the environment is Markov. The
actor has a narrower view: it receives the local feature vectors
and action masks, not the complete field and observed image.
The Markov property of the environment does not imply that each
local actor observation is itself a complete Markov state.

The initial distribution $\mu_0$ in
Eq.~\eqref{eq:mdp_definition} is understood as the distribution of the first
decision state obtained by sampling the initial image and applying the first
deterministic propagation on $\mathcal{M}^{0}$.

%%%%%%%%%%%%%%%%%%%%%%%%%%%%%%%%%%%%%%%%%%%%%%%%%%%%%%%%%%%%%%%%%%%%%%%%%%
\subsection{Adaptive action space}
\label{sec:rl_actions}
%%%%%%%%%%%%%%%%%%%%%%%%%%%%%%%%%%%%%%%%%%%%%%%%%%%%%%%%%%%%%%%%%%%%%%%%%%

Let $\mathcal{C}(\mathcal{M}^{n})$ denote the set of multiresolution
elements at which an adaptive decision can be made at step $n$.
For each candidate element
$i\in\mathcal{C}(\mathcal{M}^{n})$, the learned policy selects an action from
\begin{equation}
    \mathscr{A}_{i}
    =
    \left\{
        \mathrm{coarsen},
        \mathrm{retain},
        \mathrm{refine}
    \right\}.
    \label{eq:local_action_space}
\end{equation}

The actor samples one proposal for each candidate element. Let
$\mathscr{A}^{\mathrm{loc}}_i(\mathcal{M}^{n})$ contain the actions
left available by the local mask, and define the joint proposal
space by
\begin{equation}
    \mathscr{A}_{\mathrm{prop}}(\mathcal{M}^{n})
    =
    \prod_{i\in\mathcal{C}(\mathcal{M}^{n})}
    \mathscr{A}^{\mathrm{loc}}_i(\mathcal{M}^{n}).
    \label{eq:proposal_action_space}
\end{equation}
Thus the feasible proposal set depends on the current
representation. The symbol $\mathscr{A}$ in
Eq.~\eqref{eq:mdp_definition} denotes the collection of
these proposal vectors over admissible representations.
A sampled action is
\begin{equation}
    a^{n}
    =
    \{a_i^n\}_{i\in\mathcal{C}(\mathcal{M}^{n})}
    \in
    \mathscr{A}_{\mathrm{prop}}(\mathcal{M}^{n}).
    \label{eq:global_adaptive_action}
\end{equation}

Local eligibility does not mean that every joint proposal is
executed unchanged. For example, a proposed coarsening is
carried out only when the required sibling decisions agree.
The deterministic map $\mathcal{U}_{\mathrm{MR}}$ resolves such
proposals and produces an admissible representation. The
policy distribution is therefore a distribution over proposed
local actions; the realized tree is the result of applying
$\mathcal{U}_{\mathrm{MR}}$.

The policy is written as
\begin{equation}
    a^{n}
    \sim
    \pi_{\theta}
    \left(
        \,\cdot\,
        \middle|
        \mathcal{S}^{n}
    \right),
    \label{eq:policy_action}
\end{equation}
where $\theta$ denotes the trainable parameters. The stochastic formulation
in Eq.~\eqref{eq:policy_action} permits exploration during training; a
deterministic decision rule may subsequently be obtained from the learned
policy during evaluation.

%%%%%%%%%%%%%%%%%%%%%%%%%%%%%%%%%%%%%%%%%%%%%%%%%%%%%%%%%%%%%%%%%%%%%%%%%%
\subsection{Multiresolution update and dynamical transition}
\label{sec:rl_transition}
%%%%%%%%%%%%%%%%%%%%%%%%%%%%%%%%%%%%%%%%%%%%%%%%%%%%%%%%%%%%%%%%%%%%%%%%%%

At decision stage $n$, the field
$\widetilde{u}^{\,n+1}$ has already been propagated on
$\mathcal{M}^{n}$ and is contained in the decision state
$\mathcal{S}^{n}$.
The policy therefore acts only after this propagation has been completed.

Once an action $a^{n}$ has been selected, the propagated field and current
hierarchy are updated through
\begin{equation}
    \left(
        u^{n+1},
        \mathcal{M}^{n+1}
    \right)
    =
    \mathcal{U}_{\mathrm{MR}}
    \left(
        \widetilde{u}^{\,n+1},
        \mathcal{M}^{n},
        a^{n}
    \right).
    \label{eq:mr_update_section3}
\end{equation}
The map $\mathcal{U}_{\mathrm{MR}}$ interprets the sampled
proposals, resolves grouped coarsening, applies refinement to
surviving eligible leaves, and transfers the propagated field
to the resulting admissible representation. Several different
proposal vectors can produce the same accepted tree. The map
performs no further coefficient propagation. Because it does
not select the tree by minimizing the AMDI energy, hierarchical
admissibility here should not be confused with energy decay.

The result defines the accepted state
\begin{equation}
    \mathcal{X}^{n+1}
    =
    \left(
        u^{n+1},
        \mathcal{M}^{n+1},
        t_{n+1}
    \right).
    \label{eq:accepted_state_after_action}
\end{equation}

For all nonterminal decision stages
$n=0,\ldots,N_t-2$, the next MDP decision state is obtained by applying the
unchanged AMDI propagator to the accepted state,
\begin{equation}
    \widetilde{u}^{\,n+2}
    =
    \mathcal{D}_{\Delta t}
    \left(
        u^{n+1};
        \mathcal{M}^{n+1}
    \right),
    \label{eq:next_decision_propagation}
\end{equation}
followed by feature construction,
\begin{equation}
    \boldsymbol{\phi}^{n+1}
    =
    \mathcal{F}
    \left(
        u_{\mathrm{obs}},
        \widetilde{u}^{\,n+2},
        \mathcal{M}^{n+1},
        t_{n+2}
    \right).
    \label{eq:next_decision_features}
\end{equation}
Hence,
\begin{equation}
    \mathcal{S}^{n+1}
    =
    \left(
        u_{\mathrm{obs}},
        \widetilde{u}^{\,n+2},
        \mathcal{M}^{n+1},
        \boldsymbol{\phi}^{n+1},
        t_{n+2}
    \right),
    \qquad
    n=0,\ldots,N_t-2.
    \label{eq:deterministic_environment_transition}
\end{equation}

After the final decision at $n=N_t-1$, the adaptive update produces the
terminal accepted state
$\left(u^{N_t},\mathcal{M}^{N_t}\right)$.  No additional diffusion step is
performed.
The final adaptive decision is nevertheless included in the
terminal occupancy and switching measures. In particular, a
newly refined detail initialized to zero is counted as active
even though no subsequent coefficient update uses it.
For the MDP notation, we associate this endpoint with an
absorbing terminal state
\begin{equation}
    \mathcal{S}^{N_t}
    =
    \mathcal{S}_{\mathrm{term}}
    \left(
        \mathcal{X}^{N_t}
    \right),
    \qquad
    \mathcal{X}^{N_t}
    =
    \left(
        u^{N_t},
        \mathcal{M}^{N_t},
        t_{N_t}
    \right),
    \qquad
    V_{\psi}
    \left(
        \mathcal{S}^{N_t}
    \right)
    =
    0.
    \label{eq:terminal_mdp_state}
\end{equation}
Thus the finite horizon contains exactly $N_t$ physical propagation steps
and $N_t$ adaptive decisions, with no propagation beyond $t_{N_t}$.

The complete transition can therefore be represented by
\begin{equation}
    \mathcal{S}^{n+1}
    =
    \mathcal{G}
    \left(
        \mathcal{S}^{n},
        a^{n}
    \right),
    \qquad
    n=0,\ldots,N_t-1,
    \label{eq:environment_map_explicit}
\end{equation}
where, for $n=0,\ldots,N_t-2$, $\mathcal{G}$ consists of the constrained
adaptive update followed by the deterministic propagation and feature
construction required to form the next decision state.  For the final
decision $n=N_t-1$, $\mathcal{G}$ terminates after the adaptive update and
maps the resulting accepted state to $\mathcal{S}_{\mathrm{term}}$ without
performing any additional propagation.

For deterministic AMDI dynamics and a deterministic
multiresolution-update operator, the transition kernel is
\begin{equation}
    \mathcal{P}
    \left(
        d\mathcal{S}'
        \middle|
        \mathcal{S}^{n},
        a^{n}
    \right)
    =
    \delta_{
        \mathcal{G}(\mathcal{S}^{n},a^{n})
    }
    \left(
        d\mathcal{S}'
    \right).
    \label{eq:transition_kernel}
\end{equation}

To make the physical ordering explicit, one complete AMDI outer iteration is
\begin{equation}
    \boxed{
    \begin{aligned}
        \left(
            u^{n},
            \mathcal{M}^{n}
        \right)
        &\xrightarrow{\;\mathcal{D}_{\Delta t}\;}
        \widetilde{u}^{\,n+1},
        \\[0.4em]
        \boldsymbol{\phi}^{n}
        &=
        \mathcal{F}
        \left(
            u_{\mathrm{obs}},
            \widetilde{u}^{\,n+1},
            \mathcal{M}^{n},
            t_{n+1}
        \right),
        \\[0.4em]
        a^{n}
        &\sim
        \pi_{\theta}
        \left(
            \,\cdot\,
            \middle|
            \mathcal{S}^{n}
        \right),
        \\[0.4em]
        \left(
            u^{n+1},
            \mathcal{M}^{n+1}
        \right)
        &=
        \mathcal{U}_{\mathrm{MR}}
        \left(
            \widetilde{u}^{\,n+1},
            \mathcal{M}^{n},
            a^{n}
        \right).
    \end{aligned}
    }
    \label{eq:learned_amdi_cycle}
\end{equation}

Equation~\eqref{eq:learned_amdi_cycle} has exactly the same
propagate--then--adapt ordering as deterministic AMDI.  The only
difference is the mechanism used to select the post-propagation adaptive
update.

%%%%%%%%%%%%%%%%%%%%%%%%%%%%%%%%%%%%%%%%%%%%%%%%%%%%%%%%%%%%%%%%%%%%%%%%%%
\subsection{Accuracy, representation occupancy, and adaptive switching}
\label{sec:rl_cost_components}
%%%%%%%%%%%%%%%%%%%%%%%%%%%%%%%%%%%%%%%%%%%%%%%%%%%%%%%%%%%%%%%%%%%%%%%%%%

The learned policy should balance approximation accuracy against the size of
the active multiresolution representation. A policy that always activates
the finest-level representation can approach the uniformly resolved solution
but does not constitute a useful adaptive algorithm. Conversely, aggressive
coarsening can reduce representation occupancy at the expense of unacceptable
diffusion error. We therefore include representation occupancy explicitly in
the learning objective. Representation occupancy and adaptive switching are used here as
algorithmic complexity measures; they are not interpreted as direct
measurements of wall-clock computational cost.

Let $u_{\mathrm{ref}}^{n}$ denote the uniformly resolved reference solution
at step $n$. Because the adaptive and reference solutions may be represented
on different hierarchies, we denote by
$\mathcal{R}_{\mathrm{full}}$ the reconstruction of an adaptive field onto
the uniformly resolved comparison space and define
\begin{equation}
    \widehat{u}^{n}
    =
    \mathcal{R}_{\mathrm{full}}
    \left(
        u^{n},
        \mathcal{M}^{n}
    \right).
    \label{eq:full_reconstruction}
\end{equation}
A normalized instantaneous error is then
\begin{equation}
    E^{n}
    =
    \frac{
        \left\|
            \widehat{u}^{n}
            -
            u_{\mathrm{ref}}^{n}
        \right\|_{2}
    }{
        \left\|
            u_{\mathrm{ref}}^{n}
        \right\|_{2}
        +
        \epsilon_{\mathrm{num}}
    },
    \label{eq:instantaneous_error}
\end{equation}
where $\epsilon_{\mathrm{num}}>0$ is a numerical regularization constant.

Representation usage is measured by the occupancy
$\rho_{\mathrm{act}}^{n}$ defined in
Eq.~\eqref{eq:representation_occupancy}. This dimensionless quantity
measures the fraction of the uniformly resolved representation that remains
active at diffusion step $n$.

A third quantity is introduced to quantify changes in the active
multiresolution hierarchy between consecutive diffusion steps. Since
$\mathcal{M}^{n}$ denotes the complete adaptive representation rather than
a mathematical set, we first define the corresponding active index set by
\begin{equation}
    \Lambda^{n}
    =
    \Lambda
    \left(
        \mathcal{M}^{n}
    \right),
    \label{eq:active_index_set}
\end{equation}
where $\Lambda^{n}$ contains the indices of the multiresolution degrees of
freedom that are active at step $n$.

We then define the normalized switching measure as
\begin{equation}
    \sigma_{\mathrm{sw}}^{n}
    =
    \frac{
        \left|
            \Lambda^{n+1}
            \triangle
            \Lambda^{n}
        \right|
    }{
        N_{\mathrm{full}}
    },
    \label{eq:adaptive_switching_measure}
\end{equation}
where $\triangle$ denotes the symmetric difference between the active index
sets at two consecutive diffusion steps. Thus,
$|\Lambda^{n+1}\triangle\Lambda^{n}|$ counts the multiresolution degrees of
freedom whose activity status changes between steps $n$ and $n+1$.

The quantity $\sigma_{\mathrm{sw}}^{n}$ therefore measures the fraction of
the full representation that is activated or deactivated during one adaptive
update. The switching contribution penalizes changes in the active hierarchy within the trajectory objective.  Because the policy is optimized jointly with the accuracy and occupancy terms, this penalty need not translate into a
monotonic reduction of the realized switching measure after retraining.

%%%%%%%%%%%%%%%%%%%%%%%%%%%%%%%%%%%%%%%%%%%%%%%%%%%%%%%%%%%%%%%%%%%%%%%%%%
\subsection{Reward and trajectory objective}
\label{sec:rl_reward}
%%%%%%%%%%%%%%%%%%%%%%%%%%%%%%%%%%%%%%%%%%%%%%%%%%%%%%%%%%%%%%%%%%%%%%%%%%

Because the adaptive action deterministically generates the accepted state
$\mathcal{X}^{n+1}$ before construction of the next decision state, the
stage reward can be written as a function of the current decision state and
action,
\begin{equation}
    R:
    \mathscr{S}
    \times
    \mathscr{A}
    \longrightarrow
    \mathbb{R}.
    \label{eq:reward_map}
\end{equation}
For the action $a^{n}$ taken from $\mathcal{S}^{n}$, let
\begin{equation}
    \mathcal{X}^{n+1}
    =
    \left(
        u^{n+1},
        \mathcal{M}^{n+1},
        t_{n+1}
    \right)
\end{equation}
denote the accepted post-adaptation state generated by
$\mathcal{U}_{\mathrm{MR}}$.  The stage reward is then
\begin{equation}
    R^{n}
    =
    R
    \left(
        \mathcal{S}^{n},
        a^{n}
    \right),
\end{equation}
and is chosen as
\begin{equation}
    R^{n}
    =
    -
    \left[
        \lambda_{\mathrm{err}} E^{n+1}
        +
        \lambda_{\mathrm{occ}}
        \rho_{\mathrm{act}}^{n+1}
        +
        \lambda_{\mathrm{sw}}
        \sigma_{\mathrm{sw}}^{n}
    \right],
    \label{eq:stage_reward}
\end{equation}
where
$\lambda_{\mathrm{err}}\geq0$,
$\lambda_{\mathrm{occ}}\geq0$, and
$\lambda_{\mathrm{sw}}\geq0$,
with at least one coefficient strictly positive, control the relative
importance of reference fidelity, representation occupancy, and adaptive
switching, respectively.

For a given observed image, initial condition, and numerical
protocol, the uniformly resolved reference trajectory is
computed deterministically. The reference-dependent reward is
therefore defined for the full environment state even though
the reference field is not supplied to the actor or critic as
an input. It is used during training and evaluation, but a
deployed policy needs only the features in
Section~\ref{sec:policy_features}.

The error $E^{n+1}$ measures agreement with the uniformly
resolved AMDI computation. It does not measure agreement with
the clean image. Reconstruction quality is assessed
separately in the numerical section.

For a trajectory of $N_t$ diffusion steps,
\begin{equation}
    \tau
    =
    \left(
        \mathcal{S}^{0},
        a^{0},
        \mathcal{S}^{1},
        a^{1},
        \ldots,
        \mathcal{S}^{N_t}
    \right),
    \label{eq:trajectory_definition}
\end{equation}
the discounted return is
\begin{equation}
    G(\tau)
    =
    \sum_{n=0}^{N_t-1}
    \gamma^{n}
    R^{n}.
    \label{eq:discounted_return}
\end{equation}

The probability law of a trajectory under policy $\pi_{\theta}$ is induced
jointly by the initial-state distribution, the policy, and the transition
kernel. We denote the corresponding trajectory measure by
\begin{equation}
    \mathbb{P}_{\theta}
    \left(
        d\tau
    \right)
    =
    \mu_{0}
    \left(
        d\mathcal{S}^{0}
    \right)
    \prod_{n=0}^{N_t-1}
    \pi_{\theta}
    \left(
        da^{n}
        \middle|
        \mathcal{S}^{n}
    \right)
    \mathcal{P}
    \left(
        d\mathcal{S}^{n+1}
        \middle|
        \mathcal{S}^{n},
        a^{n}
    \right).
    \label{eq:trajectory_probability}
\end{equation}

The policy objective is the expected return
\begin{equation}
    J(\theta)
    =
    \mathbb{E}_{\tau\sim\mathbb{P}_{\theta}}
    \left[
        G(\tau)
    \right].
    \label{eq:policy_objective}
\end{equation}
PPO seeks to increase this objective from sampled trajectories;
it does not provide a global maximizer of $J$. The symbol
$\theta^\star$ used below denotes the checkpoint chosen by
validation return after training.

In the reported experiments $\gamma=1$, so rewards from all
six decision stages enter the return without discounting.
An early tree choice can affect later coefficient updates and
later rewards. This makes trajectory training meaningful as
a formulation, but it does not by itself establish that the
trained policy uses temporal feedback more effectively than
a one-step selector. That question requires a separate
numerical comparison.

%%%%%%%%%%%%%%%%%%%%%%%%%%%%%%%%%%%%%%%%%%%%%%%%%%%%%%%%%%%%%%%%%%%%%%%%%%
\subsection{Consistency with deterministic AMDI}
\label{sec:rl_consistency}
%%%%%%%%%%%%%%%%%%%%%%%%%%%%%%%%%%%%%%%%%%%%%%%%%%%%%%%%%%%%%%%%%%%%%%%%%%

The two methods have the same accepted state at the next step
if, for the same propagated field and current representation,
their adaptive updates return the same field and tree:
\begin{equation}
    \mathcal{U}_{\mathrm{MR}}
    \left(
        \widetilde{u}^{\,n+1},
        \mathcal{M}^{n},
        a^{n}
    \right)
    =
    \mathcal{U}_{\mathrm{det}}
    \left(
        \widetilde{u}^{\,n+1},
        \mathcal{M}^{n}
    \right).
    \label{eq:representation_consistency}
\end{equation}
If this equality holds at every decision and the initial
states agree, induction over the common fixed-tree update
gives identical trajectories.

This is a conditional consistency statement. It does not
assert that the factorized local policy can reproduce the
deterministic selector's global ranking of candidates, nor
that the learned tree update inherits its energy inequality.
Differences in tree selection can change the state-dependent
operator assembled at later steps.

%%%%%%%%%%%%%%%%%%%%%%%%%%%%%%%%%%%%%%%%%%%%%%%%%%%%%%%%%%%%%%%%%%%%%%%%%%
\section{Learned Policy Architecture and Optimization}
\label{sec:learning_method}
%%%%%%%%%%%%%%%%%%%%%%%%%%%%%%%%%%%%%%%%%%%%%%%%%%%%%%%%%%%%%%%%%%%%%%%%%%

The MDP formulation of Section~\ref{sec:sequential_formulation} specifies
what is learned, but it does not prescribe a particular parameterization or
optimization algorithm. We now define the implementation used for learned
AMDI. The design is guided by three requirements. First, adaptive decisions
must remain local enough to scale with the number of active multiresolution
degrees of freedom. Second, the same policy should be applicable to
representations with different numbers of active elements and, consequently,
to images at resolutions not encountered during training. Third, learning
must respect the admissibility constraints of the underlying AMDI hierarchy
rather than attempting to recover them statistically.

To satisfy these requirements, we employ a shared local policy whose
parameters are applied to every candidate multiresolution element. Locally inadmissible actions are removed by an explicit admissibility mask, and the resulting
policy is optimized with a clipped actor--critic objective. Policy-gradient and
actor--critic methods provide the underlying optimization framework
\cite{Williams1992,KondaBorkar1999}; the implementation uses proximal policy
optimization (PPO) \cite{SchulmanEtAl2017PPO}. Temporal credit
assignment is performed using generalized advantage estimation
\cite{SchulmanEtAl2016GAE}. The uniformly resolved reference trajectory enters
the training reward through Eq.~\eqref{eq:instantaneous_error}, but it is
never supplied to the policy itself.

%%%%%%%%%%%%%%%%%%%%%%%%%%%%%%%%%%%%%%%%%%%%%%%%%%%%%%%%%%%%%%%%%%%%%%%%%%
\subsection{Local feature construction}
\label{sec:policy_features}
%%%%%%%%%%%%%%%%%%%%%%%%%%%%%%%%%%%%%%%%%%%%%%%%%%%%%%%%%%%%%%%%%%%%%%%%%%

The local policy features are constructed \emph{after} the AMDI propagation
and \emph{before} the adaptive update, consistently with the propagate--then--adapt ordering of deterministic AMDI.  Starting from the accepted state
$\left(u^{n},\mathcal{M}^{n}\right)$, the field is first propagated on the
current multiresolution representation,
\begin{equation}
    \widetilde{u}^{\,n+1}
    =
    \mathcal{D}_{\Delta t}
    \left(
        u^{n};
        \mathcal{M}^{n}
    \right).
    \label{eq:feature_prepropagation}
\end{equation}
The propagated field $\widetilde{u}^{\,n+1}$ and the still-current hierarchy
$\mathcal{M}^{n}$ then define the information available to the adaptive
policy.

Let
\begin{equation}
    \mathcal{C}^{n}
    =
    \mathcal{C}
    \left(
        \mathcal{M}^{n}
    \right)
\end{equation}
denote the set of multiresolution elements at which an adaptive decision can
be made after the propagation in
Eq.~\eqref{eq:feature_prepropagation}.  For each
$i\in\mathcal{C}^{n}$, let $K_i$ denote its local spatial support and
$\ell_i$ its current resolution level.  Rather than receiving a flattened
representation of the complete image, the shared policy acts on the local
feature vector
\begin{equation}
    \mathbf{x}_{i}^{n}
    =
    \left[
        q_{i}^{n},
        g_{i}^{n},
        v_{i}^{n},
        \frac{\ell_i}{L},
        \rho_{i}^{n},
        \frac{n+1}{N_t}
    \right]^{\mathsf T},
    \label{eq:local_policy_features}
\end{equation}
where the individual components quantify the local multiresolution content,
spatial variation, current resolution level, neighborhood occupancy, and
position within the finite diffusion trajectory.

For a candidate leaf $i$, let
\begin{equation}
    \mathbf{d}_{i}^{\mathrm{obs}}
    =
    \left[
        d_{i,1}^{\mathrm{obs}},
        d_{i,2}^{\mathrm{obs}},
        d_{i,3}^{\mathrm{obs}}
    \right]^{\mathsf T}
\end{equation}
denote the three tensor-product Haar detail coefficients of the observed image
associated with the prospective refinement of that leaf.  The local
multiresolution-magnitude feature is
\begin{equation}
    q_i^n
    =
    \frac{
        \left\|
            \mathbf{d}_{i}^{\mathrm{obs}}
        \right\|_2
    }{
        s_q+\epsilon_{\mathrm{feat}}
    }.
    \label{eq:coefficient_feature}
\end{equation}
The observed-image Haar details do not evolve during an
episode; the superscript $n$ records the changing candidate
set. Changes in the policy input over time come from the
propagated-field features, hierarchy, and decision stage.
The same observed-detail magnitude also provides a natural
non-learned ranking rule for the numerical comparisons.

Let $K_i^{+}$ denote the pixel support of candidate leaf $i$ enlarged by a
one-pixel halo in each coordinate direction and clipped at the image
boundary.  The post-propagation spatial features are
\begin{equation}
    g_i^n
    =
    \frac{
        \left[
            \displaystyle
            \sum_{p\in K_i^{+}}
            \left\|
                \nabla_h
                \widetilde{u}^{\,n+1}(p)
            \right\|_2^2
        \right]^{1/2}
    }{
        s_g+\epsilon_{\mathrm{feat}}
    },
    \label{eq:gradient_feature}
\end{equation}
and
\begin{equation}
    v_i^n
    =
    \frac{
        \operatorname{Var}_{K_i^{+}}
        \left(
            \widetilde{u}^{\,n+1}
        \right)
    }{
        s_v+\epsilon_{\mathrm{feat}}
    }.
    \label{eq:variance_feature}
\end{equation}
Here $\nabla_h$ denotes the discrete image gradient used in the numerical
implementation.

The normalized hierarchy level
\begin{equation}
    \frac{\ell_i}{L}
    \in
    [0,1]
\end{equation}
provides the policy with the location of the candidate element within the
available multiresolution hierarchy, where $L$ denotes the maximum admissible
resolution level.

Let $\mathcal{L}^{n}$ denote the set of active leaf cells of the current
hierarchy.
For a candidate $i$, $\mathcal{N}_i$ contains its cell, the
geometrically present siblings, its parent when present, and
same-level face-neighbor cells inside the domain. Its
denominator counts these locations whether or not they are
active leaves; its numerator counts only members of the active
leaf set $\mathcal{L}^{n}$. In particular, an internal parent
may contribute to the denominator without contributing to
the numerator.
The local occupancy feature is
\begin{equation}
    \rho_i^n
    =
    \frac{
        \left|
            \mathcal{L}^{n}
            \cap
            \mathcal{N}_i
        \right|
    }{
        \left|
            \mathcal{N}_i
        \right|
    }.
    \label{eq:local_occupancy_feature}
\end{equation}

Finally,
\begin{equation}
    \frac{n+1}{N_t}
\end{equation}
identifies the normalized position of the post-propagation decision state
within the finite diffusion horizon.  The shift by one reflects the fact
that the policy acts after propagation from $t_n$ to $t_{n+1}$.

The complete feature component of the decision state introduced in
Eq.~\eqref{eq:rl_state_definition} is therefore
\begin{equation}
    \boldsymbol{\phi}^{n}
    =
    \left\{
        \mathbf{x}_{i}^{n}
    \right\}_{i\in\mathcal{C}^{n}}.
    \label{eq:complete_feature_state}
\end{equation}
The scales $s_q$, $s_g$, and $s_v$ are fixed before policy training as the
$0.95$ quantiles of the absolute raw values of the corresponding features
computed over the training images.  The scale-estimation trajectories use
the first three post-propagation decision states generated with deterministic
AMDI adaptation.  The same scales are subsequently used for training,
validation, testing, resolution transfer, and deployment, with
$\epsilon_{\mathrm{feat}}=10^{-12}$.  The uniformly resolved reference
trajectory is not used in the feature construction and is therefore
unavailable to the policy at inference time.

%%%%%%%%%%%%%%%%%%%%%%%%%%%%%%%%%%%%%%%%%%%%%%%%%%%%%%%%%%%%%%%%%%%%%%%%%%
\subsection{Shared local policy}
\label{sec:policy_network}
%%%%%%%%%%%%%%%%%%%%%%%%%%%%%%%%%%%%%%%%%%%%%%%%%%%%%%%%%%%%%%%%%%%%%%%%%%

The number of candidate elements varies during the diffusion trajectory and
changes with image resolution. A fixed-size policy acting on the complete
multiresolution state would therefore introduce an unnecessary dependence on
the discretization size. We instead use a shared local policy: the same
parameterized map is applied independently to the feature vector of every
candidate element.

Let
\begin{equation}
    \mathbf{h}_{i}^{n}
    =
    f_{\theta}^{\mathrm{enc}}
    \left(
        \mathbf{x}_{i}^{n}
    \right)
\end{equation}
denote the latent representation obtained from a multilayer perceptron with
shared parameters. In the implementation considered here, the encoder
contains two fully connected hidden layers,
\begin{equation}
    \mathbf{h}_{i}^{n}
    =
    \sigma
    \left[
        W_{2}
        \,
        \sigma
        \left(
            W_{1}\mathbf{x}_{i}^{n}
            +
            \mathbf{b}_{1}
        \right)
        +
        \mathbf{b}_{2}
    \right],
    \label{eq:policy_encoder}
\end{equation}
where $\sigma=\tanh$ in all numerical experiments reported below.

The actor head maps this latent representation to three logits,
\begin{equation}
    \mathbf{z}_{i}^{n}
    =
    W_{\pi}
    \mathbf{h}_{i}^{n}
    +
    \mathbf{b}_{\pi}
    \in
    \mathbb{R}^{3},
    \label{eq:policy_logits}
\end{equation}
corresponding to the actions
\begin{equation}
    \left\{
        \mathrm{coarsen},
        \mathrm{retain},
        \mathrm{refine}
    \right\}.
\end{equation}

Weight sharing across candidate elements is important for two reasons. It
makes the number of trainable parameters independent of the number of active
multiresolution degrees of freedom, and it allows the same local decision
rule to be applied to previously unseen spatial resolutions. Resolution
dependence enters explicitly through the normalized level feature
$\ell_i/L$ rather than through the dimensions of the network.

%%%%%%%%%%%%%%%%%%%%%%%%%%%%%%%%%%%%%%%%%%%%%%%%%%%%%%%%%%%%%%%%%%%%%%%%%%
\subsection{Admissibility masking}
\label{sec:action_masking}
%%%%%%%%%%%%%%%%%%%%%%%%%%%%%%%%%%%%%%%%%%%%%%%%%%%%%%%%%%%%%%%%%%%%%%%%%%

The learned policy is not permitted to generate an invalid
multiresolution hierarchy. We therefore use an action mask to remove locally
inadmissible decisions before sampling. Global consistency of the resulting
multiresolution hierarchy remains enforced deterministically by
$\mathcal{U}_{\mathrm{MR}}$, using the admissibility construction inherited
from AMDI.

For candidate element $i$, let
\begin{equation}
    m_{i,a}^{n}
    \in
    \left\{
        0,1
    \right\},
    \qquad
    a\in\mathscr{A}_{i},
    \label{eq:action_mask}
\end{equation}
where $m_{i,a}^{n}=1$ if action $a$ is admissible for the current
multiresolution representation and $m_{i,a}^{n}=0$ otherwise. The
retain action is always admissible, ensuring that at least one valid action
exists for every candidate element.

The masked categorical policy is defined by
\begin{equation}
    \pi_{\theta}
    \left(
        a_i^n=a
        \middle|
        \mathbf{x}_{i}^{n},
        \mathbf{m}_{i}^{n}
    \right)
    =
    \frac{
        m_{i,a}^{n}
        \exp
        \left(
            z_{i,a}^{n}
        \right)
    }{
        \displaystyle
        \sum_{b\in\mathscr{A}_{i}}
        m_{i,b}^{n}
        \exp
        \left(
            z_{i,b}^{n}
        \right)
    }.
    \label{eq:masked_policy}
\end{equation}

The mask removes locally ineligible proposals before sampling;
retain remains available at every candidate. The sampled
proposals are then passed to $\mathcal{U}_{\mathrm{MR}}$.
A parent is coarsened only if all four active leaf children
propose coarsening. Proposals that do not form such a group
produce no coarsening. Refinement is applied afterward to
eligible leaves that remain active. The update map, rather
than the independent local samples, guarantees that the
resulting hierarchy is admissible.

Conditioned on the complete state, the joint adaptive policy is represented
by the factorized distribution
\begin{equation}
    \pi_{\theta}
    \left(
        a^{n}
        \middle|
        \mathcal{S}^{n}
    \right)
    =
    \prod_{i\in\mathcal{C}^{n}}
    \pi_{\theta}
    \left(
        a_i^n
        \middle|
        \mathbf{x}_{i}^{n},
        \mathbf{m}_{i}^{n}
    \right).
    \label{eq:factorized_policy}
\end{equation}
Equation~\eqref{eq:factorized_policy} describes the
distribution of sampled proposals. The deterministic map
$\mathcal{U}_{\mathrm{MR}}$ can send different proposals to
the same realized hierarchy. Rewards are assigned after this
map has been applied, while policy log-probabilities are
computed for the proposals that were actually sampled.

%%%%%%%%%%%%%%%%%%%%%%%%%%%%%%%%%%%%%%%%%%%%%%%%%%%%%%%%%%%%%%%%%%%%%%%%%%
\subsection{Value function}
\label{sec:value_function}
%%%%%%%%%%%%%%%%%%%%%%%%%%%%%%%%%%%%%%%%%%%%%%%%%%%%%%%%%%%%%%%%%%%%%%%%%%

PPO is implemented in actor--critic form. Because the number of candidate
elements varies with $\mathcal{M}^{n}$, the critic must also accept a
variable-size adaptive representation. We construct a fixed-dimensional
global descriptor through permutation-invariant pooling of the local
features,
\begin{equation}
    \overline{\mathbf{x}}^{\,n}
    =
    \frac{1}{
        |\mathcal{C}^{n}|
    }
    \sum_{i\in\mathcal{C}^{n}}
    \mathbf{x}_{i}^{n},
    \label{eq:mean_feature_pool}
\end{equation}
and
\begin{equation}
    \mathbf{x}_{\max}^{\,n}
    =
    \operatorname*{max}_{i\in\mathcal{C}^{n}}
    \mathbf{x}_{i}^{n},
    \label{eq:max_feature_pool}
\end{equation}
where the maximum is taken componentwise.

The critic input is
\begin{equation}
    \mathbf{g}^{n}
    =
    \left[
        \overline{\mathbf{x}}^{\,n},
        \mathbf{x}_{\max}^{\,n},
        \rho_{\mathrm{act}}^{n},
        \frac{n+1}{N_t}
    \right].
    \label{eq:critic_global_features}
\end{equation}
The value function is represented by
\begin{equation}
    V_{\psi}
    \left(
        \mathcal{S}^{n}
    \right)
    =
    f_{\psi}^{V}
    \left(
        \mathbf{g}^{n}
    \right),
    \label{eq:value_network}
\end{equation}
where $\psi$ denotes the critic parameters. Mean and maximum pooling make the
critic independent of the cardinality of the current candidate set.

%%%%%%%%%%%%%%%%%%%%%%%%%%%%%%%%%%%%%%%%%%%%%%%%%%%%%%%%%%%%%%%%%%%%%%%%%%
\subsection{Advantage estimation}
\label{sec:advantage_estimation}
%%%%%%%%%%%%%%%%%%%%%%%%%%%%%%%%%%%%%%%%%%%%%%%%%%%%%%%%%%%%%%%%%%%%%%%%%%

For a rollout generated by the current policy, the one-step temporal
difference residual is evaluated using the critic parameters fixed during
rollout generation,
\begin{equation}
    \delta_{n}
    =
    R^{n}
    +
    \gamma
    V_{\psi_{\mathrm{old}}}
    \left(
        \mathcal{S}^{n+1}
    \right)
    -
    V_{\psi_{\mathrm{old}}}
    \left(
        \mathcal{S}^{n}
    \right).
    \label{eq:td_residual}
\end{equation}
At the terminal state,
\begin{equation}
    V_{\psi_{\mathrm{old}}}
    \left(
        \mathcal{S}^{N_t}
    \right)
    =
    V_{\psi_{\mathrm{old}}}
    \left(
        \mathcal{S}_{\mathrm{term}}
    \right)
    =
    0.
\end{equation}

Generalized advantage estimation is then used to construct
\begin{equation}
    \widehat{A}_{n}
    =
    \sum_{l=0}^{N_t-n-1}
    \left(
        \gamma\lambda_{\mathrm{GAE}}
    \right)^{l}
    \delta_{n+l},
    \label{eq:gae_advantage}
\end{equation}
where
$\lambda_{\mathrm{GAE}}\in[0,1]$ controls the bias--variance tradeoff
\cite{SchulmanEtAl2016GAE}. The corresponding value-function target is
\begin{equation}
    \widehat{V}_{n}
    =
    \widehat{A}_{n}
    +
    V_{\psi_{\mathrm{old}}}
    \left(
        \mathcal{S}^{n}
    \right),
    \label{eq:value_target}
\end{equation}
where $\psi_{\mathrm{old}}$ denotes the critic parameters used to generate
the rollout.

%%%%%%%%%%%%%%%%%%%%%%%%%%%%%%%%%%%%%%%%%%%%%%%%%%%%%%%%%%%%%%%%%%%%%%%%%%
\subsection{Clipped policy optimization}
\label{sec:ppo_optimization}
%%%%%%%%%%%%%%%%%%%%%%%%%%%%%%%%%%%%%%%%%%%%%%%%%%%%%%%%%%%%%%%%%%%%%%%%%%

Because the adaptive action $a^{n}$ consists of the complete set of local
decisions taken at diffusion step $n$, the PPO likelihood ratio must be
defined for the corresponding joint action. Using the factorization in
Eq.~\eqref{eq:factorized_policy}, the log-probability ratio is
\begin{equation}
    \log r_{n}(\theta)
    =
    \sum_{i\in\mathcal{C}^{n}}
    \left[
        \log
        \pi_{\theta}
        \left(
            a_{i}^{n}
            \middle|
            \mathbf{x}_{i}^{n},
            \mathbf{m}_{i}^{n}
        \right)
        -
        \log
        \pi_{\theta_{\mathrm{old}}}
        \left(
            a_{i}^{n}
            \middle|
            \mathbf{x}_{i}^{n},
            \mathbf{m}_{i}^{n}
        \right)
    \right].
    \label{eq:ppo_log_probability_ratio}
\end{equation}
The probability ratio for the complete sampled proposal is
\begin{equation}
    r_n(\theta)
    =
    \exp\!\left[\log r_n(\theta)\right].
    \label{eq:ppo_probability_ratio}
\end{equation}
This ratio concerns the sampled proposal vector, not the tree
produced by $\mathcal{U}_{\mathrm{MR}}$. Its log-ratio is summed
over candidates without dividing by $|\mathcal{C}^{n}|$.

For the actor update, let $\widetilde{A}_n$ denote
$\widehat{A}_n$ standardized over the decision-step samples
in the pooled rollout. The value target in
Eq.~\eqref{eq:value_target} retains the unstandardized
advantage. The clipped PPO objective is
\begin{equation}
    \mathcal{L}_{\mathrm{clip}}(\theta)
    =
    \mathbb{E}\!\left[
        \min\!\left(
            r_n(\theta)\widetilde{A}_n,\,
            \operatorname{clip}\!\left(
                r_n(\theta),
                1-\epsilon_{\mathrm{PPO}},
                1+\epsilon_{\mathrm{PPO}}
            \right)\widetilde{A}_n
        \right)
    \right].
    \label{eq:ppo_clip_objective}
\end{equation}
Here $\epsilon_{\mathrm{PPO}}>0$ specifies the clipping
interval for a joint decision step
\cite{SchulmanEtAl2017PPO}. Since the number of candidates
can change between steps, a fixed joint clipping interval
does not impose the same bound on each local action.

The critic is optimized using
\begin{equation}
    \mathcal{L}_{V}(\psi)
    =
    \mathbb{E}
    \left[
        \left(
            V_{\psi}
            \left(
                \mathcal{S}^{n}
            \right)
            -
            \widehat{V}_{n}
        \right)^{2}
    \right].
    \label{eq:value_loss}
\end{equation}

To maintain sufficient exploration during training, we also define the
mean categorical entropy
\begin{equation}
    \mathcal{H}(\theta)
    =
    \mathbb{E}
    \left[
        \frac{1}{
            |\mathcal{C}^{n}|
        }
        \sum_{i\in\mathcal{C}^{n}}
        H
        \left[
            \pi_{\theta}
            \left(
                \cdot
                \middle|
                \mathbf{x}_{i}^{n},
                \mathbf{m}_{i}^{n}
            \right)
        \right]
    \right].
    \label{eq:policy_entropy}
\end{equation}

The complete objective minimized during training is
\begin{equation}
    \mathcal{L}_{\mathrm{train}}
    =
    -
    \mathcal{L}_{\mathrm{clip}}
    +
    c_{V}
    \mathcal{L}_{V}
    -
    c_{H}
    \mathcal{H},
    \label{eq:complete_training_loss}
\end{equation}
where $c_V\geq0$ and $c_H\geq0$ control the contributions of the critic loss
and entropy regularization, respectively.
Optimization uses Adam~\cite{KingmaBa2015}. A PPO update in
the reported six-step protocol collects 16 trajectories,
giving 96 decision-step samples. In each of four optimization
epochs, the steps are shuffled and processed with a
minibatch of one decision step, with gradient clipping at
each optimizer step. Thus one PPO update comprises 384 Adam
steps under the reported protocol.

%%%%%%%%%%%%%%%%%%%%%%%%%%%%%%%%%%%%%%%%%%%%%%%%%%%%%%%%%%%%%%%%%%%%%%%%%%
\subsection{Training protocol}
\label{sec:training_protocol}
%%%%%%%%%%%%%%%%%%%%%%%%%%%%%%%%%%%%%%%%%%%%%%%%%%%%%%%%%%%%%%%%%%%%%%%%%%

Training episodes are generated by sampling initial images from the prescribed
training-image distribution.  After construction of the initial
multiresolution representation $\mathcal{M}^{0}$ and the first deterministic
propagation, this image distribution induces the initial decision-state
distribution $\mu_0$ appearing in Eq.~\eqref{eq:mdp_definition}. For every training initial condition, the
uniformly resolved diffusion trajectory required to evaluate
$E^n$ is computed independently using the same diffusion parameters and
temporal discretization as the adaptive calculation. These reference
trajectories are used only to construct the reward and are not supplied to
the actor or critic as state information.

One training episode proceeds as follows:
\begin{enumerate}

    \item Sample an initial image and construct the initial accepted
    representation $\mathcal{M}^{0}$ and field $u^{0}$.

    \item For each outer iteration $n$, first propagate the accepted field on
    the current representation,
    \begin{equation}
        \widetilde{u}^{\,n+1}
        =
        \mathcal{D}_{\Delta t}
        \left(
            u^{n};
            \mathcal{M}^{n}
        \right).
    \end{equation}

    \item Construct the post-propagation decision state
    $\mathcal{S}^{n}$, the candidate set $\mathcal{C}^{n}$, the local
    feature vectors $\mathbf{x}_{i}^{n}$, and the admissibility masks
    $\mathbf{m}_{i}^{n}$.

    \item Sample the adaptive actions from the masked policy
    $\pi_{\theta}$.

    \item Apply the selected actions to the propagated state using
    $\mathcal{U}_{\mathrm{MR}}$,
    \begin{equation}
        \left(
            u^{n+1},
            \mathcal{M}^{n+1}
        \right)
        =
        \mathcal{U}_{\mathrm{MR}}
        \left(
            \widetilde{u}^{\,n+1},
            \mathcal{M}^{n},
            a^{n}
        \right).
    \end{equation}

    \item Evaluate the approximation error, representation occupancy, and
    switching measure on the resulting accepted state and construct the
    reward using Eq.~\eqref{eq:stage_reward}.

    \item Continue until $n=N_t-1$, producing a complete rollout.

    \item Compute the advantages using
    Eq.~\eqref{eq:gae_advantage} and update the actor and critic using
    Eq.~\eqref{eq:complete_training_loss}.

\end{enumerate}

In the reported protocol, $\lambda_{\mathrm{err}}=1$;
$\lambda_{\mathrm{occ}}$ and $\lambda_{\mathrm{sw}}$ set the
relative penalties for representation occupancy and tree
changes. The architecture and optimizer settings are stated
in Section~\ref{sec:numerical_experiments}. For each training
run, we select the checkpoint with the highest mean return
under greedy actions on the validation images. The held-out
test images are used only for final evaluation.

Because $E^n$, $\rho_{\mathrm{act}}^n$, and
$\sigma_{\mathrm{sw}}^n$ are dimensionless, the reward coefficients have a
direct interpretation as relative tradeoff parameters rather than compensating
for incompatible physical units.

%%%%%%%%%%%%%%%%%%%%%%%%%%%%%%%%%%%%%%%%%%%%%%%%%%%%%%%%%%%%%%%%%%%%%%%%%%
\subsection{Deployment and resolution transfer}
\label{sec:policy_deployment}
%%%%%%%%%%%%%%%%%%%%%%%%%%%%%%%%%%%%%%%%%%%%%%%%%%%%%%%%%%%%%%%%%%%%%%%%%%

During deployment, the propagate--then--adapt ordering of deterministic AMDI~\cite{TantardiniAMDI2026} is preserved.
At each outer iteration, the field is first propagated on the current
representation,
\begin{equation}
    \widetilde{u}^{\,n+1}
    =
    \mathcal{D}_{\Delta t}
    \left(
        u^{n};
        \mathcal{M}^{n}
    \right).
    \label{eq:deployment_propagation}
\end{equation}
The local features and admissibility masks are then constructed from the
post-propagation state.  Stochastic exploration is removed, and for each
candidate element the action is selected according to
\begin{equation}
    a_{i}^{n}
    =
    \arg\max_{
        a\in\mathscr{A}_{i}
    }
    \pi_{\theta^{\star}}
    \left(
        a
        \middle|
        \mathbf{x}_{i}^{n},
        \mathbf{m}_{i}^{n}
    \right).
    \label{eq:deterministic_deployment_policy}
\end{equation}
The selected actions are subsequently applied through
$\mathcal{U}_{\mathrm{MR}}$ to obtain
$\left(u^{n+1},\mathcal{M}^{n+1}\right)$.

No uniformly resolved reference trajectory is required during deployment,
and the critic is not needed once training has been completed.  Deployment
therefore requires only the unchanged AMDI propagation, deterministic feature
construction, actor evaluation, admissibility masking, and the constrained
post-propagation multiresolution update.

The shared local parameterization makes the number of actor parameters
independent of both $N_{\mathrm{act}}^{n}$ and $N_{\mathrm{full}}$. Together
with the normalized level and time features, this permits the same policy to
be evaluated on multiresolution representations containing different numbers
of active elements. This architectural property does not by itself guarantee
generalization to unseen resolutions; rather, it removes a fixed-dimensional
restriction that would otherwise prevent such transfer.
Resolution transfer is therefore treated as an empirical generalization
assessment and is evaluated explicitly in
Section~\ref{sec:generalization}.

The operational order during deployment is
\begin{equation}
    \boxed{
    \mathcal{D}_{\Delta t}
    \;\longrightarrow\;
    \mathcal{F}
    \;\longrightarrow\;
    \pi_{\theta^\star}
    \;\longrightarrow\;
    \mathcal{U}_{\mathrm{MR}}
    }.
    \label{eq:learned_deterministic_separation}
\end{equation}
The field is updated by the specified fixed-tree AMDI
procedure, features and masks are constructed
deterministically, and the actor proposes local tree
actions. The final map enforces hierarchy admissibility.
The critic and uniformly resolved reference are not needed
during deployment. Because the operator assembled for a
later step depends on the accepted tree, retaining the AMDI
update procedure does not mean that two different adaptive
trajectories use identical realized operators.

%%%%%%%%%%%%%%%%%%%%%%%%%%%%%%%%%%%%%%%%%%%%%%%%%%%%%%%%%%%%%%%%%%%%%%%%%%
\section{Numerical validation}
\label{sec:numerical_experiments}
%%%%%%%%%%%%%%%%%%%%%%%%%%%%%%%%%%%%%%%%%%%%%%%%%%%%%%%%%%%%%%%%%%%%%%%%%%

We assess Learned AMDI against deterministic AMDI, a uniformly resolved AMDI
reference, and nonlearned controls constructed from the observed-image Haar
details.  The numerical study addresses six questions: (i) whether the
modular implementation reproduces the retained deterministic implementation;
(ii) which actions the deterministic and learned selectors propose and
execute; (iii) whether a learned selector improves on validation-tuned
one-step detail rules at comparable occupancy; (iv) how the
accuracy--representation tradeoff changes with the occupancy penalty;
(v) whether the selected policy and the frozen controls transfer to finer
resolutions; and (vi) how sensitive the conclusions are to policy
initialization and to the reward terms.

All methods use the same AMDI propagator, Haar representation, frozen-weight
update, fixed-tree energy safeguard, and hierarchical admissibility
constraints.  They differ only in the post-propagation treatment of the
adaptive tree.  Deterministic AMDI~\cite{TantardiniAMDI2026} minimizes its prescribed energy-plus-tree
penalty over retain, refinement, and coarsening candidates.  Learned AMDI uses
the frozen PPO actor defined in Section~\ref{sec:learning_method}.  The
additional controls either retain the initial tree or refine once according
to the observed-image Haar details and then retain the resulting tree.

%%%%%%%%%%%%%%%%%%%%%%%%%%%%%%%%%%%%%%%%%%%%%%%%%%%%%%%%%%%%%%%%%%%%%%%%%%
\subsection{Synthetic data and evaluation sets}
\label{sec:synthetic_dataset}
%%%%%%%%%%%%%%%%%%%%%%%%%%%%%%%%%%%%%%%%%%%%%%%%%%%%%%%%%%%%%%%%%%%%%%%%%%

The experiments use a deterministic, seed-controlled family of multiscale
images,
\begin{equation}
    u_{\mathrm{clean}}
    =
    \mathcal{G}_{\mathrm{RMS}}(N;s),
    \label{eq:random_multiscale_generator}
\end{equation}

where \(N\times N\) is the image resolution and \(s\) initializes the random
number generator. At the pixel centers
\[
X_{ij}=\frac{j+1/2}{N},
\qquad
Y_{ij}=\frac{i+1/2}{N},
\]
the unnormalized image is initialized as
\[
g(X,Y)=0.15+0.20X+0.10Y.
\]
Three Gaussian components
\[
A_k
\exp\left[
-\frac{(X-x_k)^2+(Y-y_k)^2}{q_k}
\right],
\qquad k=1,2,3,
\]
are then added, with
\[
(x_k,y_k)\sim\mathcal{U}([0.15,0.85]^2),
\qquad
q_k\sim\mathcal{U}(0.025,0.12),
\qquad
A_k\sim\mathcal{U}(0.15,0.45).
\]
A rectangular inclusion has lower-left corner sampled from
\(\mathcal{U}([0.1,0.6]^2)\), independent side lengths sampled from
\(\mathcal{U}(0.12,0.30)\), and amplitude sampled from
\(\mathcal{U}(0.15,0.35)\). Finally, a local oscillatory contribution
\[
0.10\sin(2\pi fX)\sin(2\pi fY)
\]
is added on the square
\(\lvert X-p_x\rvert<0.20\), \(\lvert Y-p_y\rvert<0.20\), where
\((p_x,p_y)\sim\mathcal{U}([0.45,0.75]^2)\) and
\(f\) is sampled uniformly from \(\{4,\ldots,9\}\). The result is normalized
by
\[
u_{\mathrm{clean}}
=
\frac{g-\min g}{\max g-\min g}.
\]
All draws are generated deterministically from \(s\), so a fixed pair
\((N,s)\) reproduces the same image.

The independently prescribed four-region target is
\begin{equation}
u_{\mathrm{4reg}}(X,Y)
=
\begin{cases}
0.25,
& X<0.5,\ Y<0.5,\\
0.15+1.4(X-0.5),
& X\geq0.5,\ Y<0.5,\\
0.90,
& (X-0.25)^2+(Y-0.75)^2\leq0.12^2,\\
0.15,
& X<0.5,\ Y\geq0.5
  \text{ outside the disk},\\
0.50+0.28\sin(18\pi X)\sin(18\pi Y),
& X\geq0.5,\ Y\geq0.5.
\end{cases}
\label{eq:four_region_target}
\end{equation}

The observed image is
\begin{equation}
    u_{\mathrm{obs}}
    =
    \Pi_{[0,1]}
    \left(u_{\mathrm{clean}}+\eta\right),
    \qquad
    \eta_{ij}\overset{\mathrm{iid}}{\sim}\mathcal{N}(0,\sigma^2),
    \qquad
    \sigma=0.08,
    \label{eq:synthetic_noise_model}
\end{equation}
with pointwise clipping $\Pi_{[0,1]}$.

Training is performed only at $N_{\mathrm{train}}=32$.  The image-seed sets
are
\begin{equation}
\begin{aligned}
    \mathcal{I}_{\mathrm{train}}&=\{0,1,\ldots,15\},\\
    \mathcal{I}_{\mathrm{val}}&=\{100,101,102,103\},\\
    \mathcal{I}_{\mathrm{test}}&=\{1000,1001,\ldots,1007\}.
\end{aligned}
\label{eq:data_seed_splits}
\end{equation}
The sets are mutually disjoint.  Validation images are used for checkpoint
selection and for selecting the parameters of the nonlearned controls; test
images are not used for either purpose.

The principal holdout set contains the eight images generated from
$\mathcal{I}_{\mathrm{test}}$ and one independently prescribed four-region
benchmark with constant, smoothly varying, edge-dominated, and textured
regions.  Thus the holdout comparison contains nine cases.  The Pareto
comparison uses the common four-image subset
\begin{equation}
    \mathcal{I}_{4}=\{1000,1001,1002,1003\}
    \subset\mathcal{I}_{\mathrm{test}},
    \label{eq:four_image_subset}
\end{equation}
with the same noise realization for every method at a given image seed.

Resolution transfer is evaluated on the same four image seeds at
$N\in\{32,64,128\}$.  The $32\times32$ set is an independent-noise
same-resolution deployment check, while the $64\times64$ and
$128\times128$ sets are genuine transfer calculations.  Neither the actor nor
the top-$K$ and threshold parameters is retuned at the two finer resolutions.

Two references are used and should be distinguished.  The clean image
$u_{\mathrm{clean}}$ defines RMSE and the structural similarity index (SSIM)
\cite{WangEtAl2004SSIM}.  The uniformly resolved AMDI
trajectory defines the algorithmic error $E_{\mathrm{ref}}$ in
Eq.~\eqref{eq:instantaneous_error}.  Consequently, agreement with uniform
AMDI and agreement with the clean image are different objectives.

%%%%%%%%%%%%%%%%%%%%%%%%%%%%%%%%%%%%%%%%%%%%%%%%%%%%%%%%%%%%%%%%%%%%%%%%%%
\subsection{Numerical protocol and implementation checks}
\label{sec:numerical_protocol}
%%%%%%%%%%%%%%%%%%%%%%%%%%%%%%%%%%%%%%%%%%%%%%%%%%%%%%%%%%%%%%%%%%%%%%%%%%

Each trajectory contains six propagate--then--adapt outer iterations.  The
AMDI parameters are fixed in all comparisons:
\begin{equation}
\begin{aligned}
    \alpha &=0.010,
    &
    \beta &=2\times10^{-4},
    &
    \tau &=10^{-7},\\
    \sigma_c&=0.12,
    &
    h&=0.60,
    &
    \zeta&=10^{-6}.
\end{aligned}
\label{eq:fixed_amdi_parameters}
\end{equation}
The remaining energy and graph parameters are
\[
\gamma_{\mathrm{lev}}=0.1,
\qquad
\sigma_x=2,
\qquad
\sigma_{\ell}=1.5,
\qquad
\kappa_{\mathrm{ref}}=0.35,
\]
where a parent--child graph weight is multiplied by
\(2^{-\kappa_{\mathrm{ref}}\max(\ell_i,\ell_j)}\). Graph weights no larger
than \(10^{-14}\) are discarded, and the nonsmoothed \(\ell^1\) sparsity term
is used. The numerical and feature-normalization constants are
\(\epsilon_{\mathrm{num}}=\epsilon_{\mathrm{feat}}=10^{-12}\). The initial
Haar-detail threshold is \(0.012\), the minimum admissible level is \(2\),
and the deterministic refinement and coarsening fractions are both \(0.10\).

The actor and critic each have two width-$64$ hidden layers with hyperbolic
tangent activations.  PPO uses
\begin{equation}
\begin{aligned}
    \gamma&=1,
    &
    \lambda_{\mathrm{GAE}}&=0.95,
    &
    \epsilon_{\mathrm{PPO}}&=0.10,\\
    c_V&=0.5,
    &
    c_H&=0.02,
    &
    \eta_{\mathrm{Adam}}&=5\times10^{-5}.
\end{aligned}
\label{eq:final_ppo_parameters}
\end{equation}
Each update contains $16$ complete trajectories and four optimization epochs,
with gradient norm clipped at $0.5$.  Training is run for at most $400$
updates.  For each training seed, the checkpoint with the largest mean
validation return under deterministic action selection is retained.  The
principal checkpoint has training seed $20260811$ and reward coefficients
\begin{equation}
    \lambda_{\mathrm{err}}=1,
    \qquad
    \lambda_{\mathrm{occ}}=0.15,
    \qquad
    \lambda_{\mathrm{sw}}=0.02.
    \label{eq:final_reward_parameters}
\end{equation}

We perform two distinct deterministic regression checks.  First, a low-level
backend comparison applies the direct deterministic solver and its modular
wrapper to the same state.  Both paths return the same tree with $190$ active
degrees of freedom,
\begin{equation}
    \varepsilon_{\mathrm{reg}}^{u}=0,
    \qquad
    d_{\mathrm{tree}}=0,
    \qquad
    C_{\mathrm{rel}}=\frac{190}{1024}=0.185546875,
    \label{eq:implementation_regression}
\end{equation}
and a maximum coefficient difference no larger than $2.2\times10^{-19}$
across the retained platform checks.  Second, a within-backend trajectory
check compares the retained deterministic AMDI pathway with the modular
propagate--then--adapt control.  In that separate case,
\begin{equation}
    C_{\mathrm{rel}}^{\mathrm{base}}
    =C_{\mathrm{rel}}^{\mathrm{ctrl}}
    =0.15625,
    \qquad
    d_{\mathrm{tree}}=0,
    \label{eq:implementation_regression_crel}
\end{equation}
and the reconstructed fields and complete trajectories coincide.
The occupancies $0.18555$ and $0.15625$ therefore belong to different
regression cases and are not interchangeable.

An end-to-end integration test applies the selected checkpoint and publication
configuration to all nine holdout cases. It reproduces the archived terminal
metrics within the prescribed tolerance and records nonzero executed learned
actions. The calculation exercises the actor, feature construction, masked
action selection, constrained tree update, and terminal metrics in the
released implementation, complementing the two deterministic regression
checks.

The uniformly resolved AMDI reference uses all $1024$ Haar degrees of freedom
at $32\times32$.  Its energy is monotone and all fixed-tree safeguards are
accepted.  For an adaptive trajectory, the reported terminal quantities are
\begin{equation}
    E_{\mathrm{ref}}:=E^{N_t},
    \qquad
    C_{\mathrm{rel}}
    :=\rho_{\mathrm{act}}^{N_t}
    =\frac{N_{\mathrm{act}}^{N_t}}{N_{\mathrm{full}}},
    \label{eq:reported_terminal_metrics}
\end{equation}
and the trajectory-mean switching measure is
\begin{equation}
    \overline{\sigma}_{\mathrm{sw}}
    =\frac{1}{N_t}\sum_{n=0}^{N_t-1}\sigma_{\mathrm{sw}}^n.
    \label{eq:mean_switching_metric}
\end{equation}
Before RMSE and SSIM are evaluated, every terminal reconstruction is clipped
pointwise to $[0,1]$ and then compared with $u_{\mathrm{clean}}$.  Unless
otherwise stated, reported evaluation-set values are arithmetic means of the
case-level quantities.

%%%%%%%%%%%%%%%%%%%%%%%%%%%%%%%%%%%%%%%%%%%%%%%%%%%%%%%%%%%%%%%%%%%%%%%%%%
\subsection{Holdout controls and step-resolved action--energy analysis}
\label{sec:holdout_validation}
%%%%%%%%%%%%%%%%%%%%%%%%%%%%%%%%%%%%%%%%%%%%%%%%%%%%%%%%%%%%%%%%%%%%%%%%%%

Three nonlearned controls are evaluated with the same cases and noise as the
selected policy.  The retain-tree control never changes the
threshold-initialized tree.  The top-$K$ control ranks eligible leaves by the
raw Euclidean norm of their three observed-image Haar details and refines the
$K$ largest at decision~1.  The threshold control refines every eligible leaf
whose raw detail norm exceeds a fixed threshold, also only at decision~1.
Both detail controls retain their trees for decisions~2--6 while continuing
the same six AMDI propagation steps. Their coefficient-ranking and
coefficient-thresholding constructions are motivated by classical spatially
adaptive wavelet denoising \cite{DonohoJohnstone1994,Donoho1995}.

The control parameters are selected using only the four validation images.
The selected policy has mean validation occupancy $0.236816$.  Maximizing the
validation trajectory return gives $K=32$ and
$\tau_{\mathrm{detail}}=6.87534\times10^{-3}$, whereas matching the learned
validation occupancy gives $K=36$ and
$\tau_{\mathrm{detail}}=6.12571\times10^{-3}$.  These four choices are frozen
before the holdout, Pareto, and transfer evaluations.

Table~\ref{tab:holdout_main} reports the resulting nine-case means and the
accepted tree changes.  Deterministic AMDI accepts no refinement in any of
the $54$ holdout decisions.  It accepts four coarsened parent groups in the
first decision of the four-region case and retains the tree in the other
$53$ decisions.  Its terminal metrics are therefore nearly identical to
those of the retain-tree control.  In the present configuration,
deterministic AMDI should consequently be interpreted as an almost static
initial-tree calculation, not as a strongly adapting baseline.

\begin{table*}[t]
\caption{
Nine-case holdout comparison.  $\overline{N}_{\mathrm{act}}$ is the mean
terminal number of active degrees of freedom.  $N_{\mathrm{ref}}/N_{\mathrm{crs}}$
gives the total number of executed refined cells and coarsened parent groups,
respectively, over nine cases and six decisions.  ``Val.-return'' and
``val.-occupancy'' denote parameters selected exclusively on the validation
images.  All terminal fields are clipped to $[0,1]$ before RMSE and SSIM are
computed.
}
\label{tab:holdout_main}
\centering
\small
\setlength{\tabcolsep}{4.0pt}
\begin{tabular}{lrrrrrrr}
\hline
Method
& $E_{\mathrm{ref}}$
& $C_{\mathrm{rel}}$
& $\overline{N}_{\mathrm{act}}$
& RMSE
& SSIM
& $\overline{\sigma}_{\mathrm{sw}}$
& $N_{\mathrm{ref}}/N_{\mathrm{crs}}$
\\
\hline
Deterministic AMDI
& 0.17496 & 0.13737 & 140.67 & 0.05391 & 0.80851 & 0.00022 & $0/4$ \\
Retain initial tree
& 0.17494 & 0.13867 & 142.00 & 0.05392 & 0.80813 & 0       & $0/0$ \\
Top-$K$, val.-return ($K=32$)
& 0.14442 & 0.23242 & 238.00 & 0.05286 & 0.80829 & 0.01563 & $288/0$ \\
Top-$K$, val.-occupancy ($K=36$)
& 0.14175 & 0.24414 & 250.00 & 0.05290 & 0.80628 & 0.01758 & $324/0$ \\
Threshold, val.-return
& 0.14617 & 0.22331 & 228.67 & 0.05178 & 0.81162 & 0.01411 & $260/0$ \\
Threshold, val.-occupancy
& 0.13792 & 0.26042 & 266.67 & 0.05266 & 0.80196 & 0.02029 & $374/0$ \\
Learned, decision~1 only
& 0.13742 & 0.26335 & 269.67 & 0.05279 & 0.80135 & 0.02078 & $383/0$ \\
Learned AMDI
& 0.13657 & 0.26660 & 273.00 & 0.05306 & 0.79899 & 0.02132 & $393/0$ \\
Uniform AMDI reference
& 0       & 1       & 1024.00& 0.07141 & 0.70150 & 0       & -- \\
\hline
\end{tabular}
\end{table*}

Table~\ref{tab:holdout_actions} resolves the behavior of the two selectors.
At decision~1, the learned actor proposes and executes \(383\) refinements.
It executes ten additional refinements at decision~2 and none thereafter.
The actor later proposes coarsening actions, but no coarsening group satisfies
the execution rules. As predicted by the structural argument in
Section~\ref{sec:deterministic_amdi}, the deterministic refinement candidate
has higher AMDI energy than the unchanged candidate in all \(54\) decisions,
even before the nonnegative tree-change penalty is added. The candidate-energy
analysis therefore confirms the refinement barrier of the implemented
construction.

\begin{table}[t]
\caption{
Nine-case action and energy analysis by adaptation decision. Here R/F/C
denotes retain/refine/coarsen for deterministic AMDI, whereas R/C denotes
learned refinement/coarsening counts. Learned proposals and executions are
totals over the nine cases. ``Energy up'' reports the number of cases in which
the post-adaptation energy exceeds the post-propagation energy.
}
\label{tab:holdout_actions}
\centering
\scriptsize
\setlength{\tabcolsep}{2.5pt}
\renewcommand{\arraystretch}{1.12}
\begin{tabular}{@{}crrrrc@{}}
\hline
Decision
& AMDI R/F/C
& Proposed R/C
& Executed R/C
& $\overline{N}_{\mathrm{act}}$
& Energy up
\\
\hline
1 & $8/0/1$ & $383/0$ & $383/0$ & 269.67 & $9/9$ \\
2 & $9/0/0$ & $10/0$  & $10/0$  & 273.00 & $6/9$ \\
3 & $9/0/0$ & $0/0$   & $0/0$   & 273.00 & $0/9$ \\
4 & $9/0/0$ & $0/5$   & $0/0$   & 273.00 & $0/9$ \\
5 & $9/0/0$ & $0/10$  & $0/0$   & 273.00 & $0/9$ \\
6 & $9/0/0$ & $0/30$  & $0/0$   & 273.00 & $0/9$ \\
\hline
\end{tabular}
\end{table}

The validation-tuned threshold rule and the selected learned policy
operate at nearby but not identical occupancies.  To preserve the paired
design, Table~\ref{tab:paired_controls} reports case-level
control-minus-learned differences; ``$\pm$'' denotes the sample standard
deviation over the nine paired differences.  The threshold control uses
a slightly smaller mean representation
($\Delta C_{\mathrm{rel}}=-0.00618$) and has a mean
$E_{\mathrm{ref}}$ only $1.35\times10^{-3}$ above Learned AMDI.
It gives smaller RMSE and larger SSIM on average and improves both
metrics in eight of the nine cases.  Learned AMDI gives the smaller
$E_{\mathrm{ref}}$ in eight cases, but this sign pattern does not isolate
the quality of spatial selection because the learned policy also uses
the larger representation.  The small reference-error offset is
consistent with the accuracy--occupancy trend shown by the control
frontiers.  The experiment therefore supports PPO as a learned
allocation mechanism, but does not establish superior spatial selection
relative to the tuned observed-detail rule.

\begin{table}[t]
\caption{
Paired holdout differences relative to Learned AMDI. Entries are
control minus learned and are reported as the mean $\pm$ sample standard
deviation over nine cases. Both controls use parameters selected by
matching the learned validation occupancy.
}
\label{tab:paired_controls}
\centering
\scriptsize
\setlength{\tabcolsep}{4pt}
\renewcommand{\arraystretch}{1.15}
\begin{tabular}{lcc}
\hline
Difference
& Top-$K$
& Threshold
\\
\hline
$\Delta E_{\mathrm{ref}}$
& $+0.00518\pm0.00958$
& $+0.00135\pm0.00106$
\\
$\Delta C_{\mathrm{rel}}$
& $-0.02246\pm0.03944$
& $-0.00618\pm0.00473$
\\
$\Delta\mathrm{RMSE}$
& $-0.00017\pm0.00148$
& $-0.00041\pm0.00037$
\\
$\Delta\mathrm{SSIM}$
& $+0.00728\pm0.01323$
& $+0.00296\pm0.00297$
\\
\hline
\end{tabular}
\end{table}

Using the same six-step scalar reward employed for policy selection, the mean
holdout returns are \(-1.14552\) for Learned AMDI, \(-1.14540\) for the
occupancy-matched threshold rule, and \(-1.14678\) for the occupancy-matched
top-\(K\) rule. The validation-return settings give \(-1.14701\) for the
threshold rule and \(-1.14758\) for top-\(K\); larger values are better. Thus,
Learned AMDI and the occupancy-matched threshold rule are also essentially
indistinguishable under the reported scalar trajectory objective.

The temporal control reaches
$E_{\mathrm{ref}}=0.13742$ when the actor is allowed to act only at
decision~1, compared with $0.13657$ under unrestricted deployment.  The mean
difference is $8.49\times10^{-4}$, obtained at $3.33$ additional active
degrees of freedom for the unrestricted policy.  The decision-1-only control
also gives a slightly smaller RMSE ($0.05279$ versus $0.05306$) and larger
SSIM ($0.80135$ versus $0.79899$).  On these static test images, the learned
behavior is therefore strongly front-loaded: the first decision produces the
substantive tree change, the second provides a small correction, and
decisions~3--6 do not execute further changes.  We consequently do not infer
material sequential control from the trajectory plot alone.  Establishing a
stronger sequential advantage would require a problem whose localized
features change during the trajectory.

The energy log also distinguishes the fixed-tree propagation safeguard from
the learned tree update.  All fixed-tree safeguards are accepted, but learned
adaptation raises the post-propagation AMDI energy in all nine cases at
decision~1 and in six cases at decision~2.  The corresponding mean
post-adaptation increments are $3.29\times10^{-5}$ and
$4.78\times10^{-7}$.  At decision~1 the complete propagation-plus-adaptation
update also increases the energy relative to the state before propagation by
$1.82\times10^{-5}$ on average.  Hence the fixed-tree safeguard guarantees
the accepted propagation update only; it does not guarantee monotone decay of
the learned outer iteration after a tree change.

\begin{figure}[t]
    \centering
    \includegraphics[width=0.98\textwidth]
    {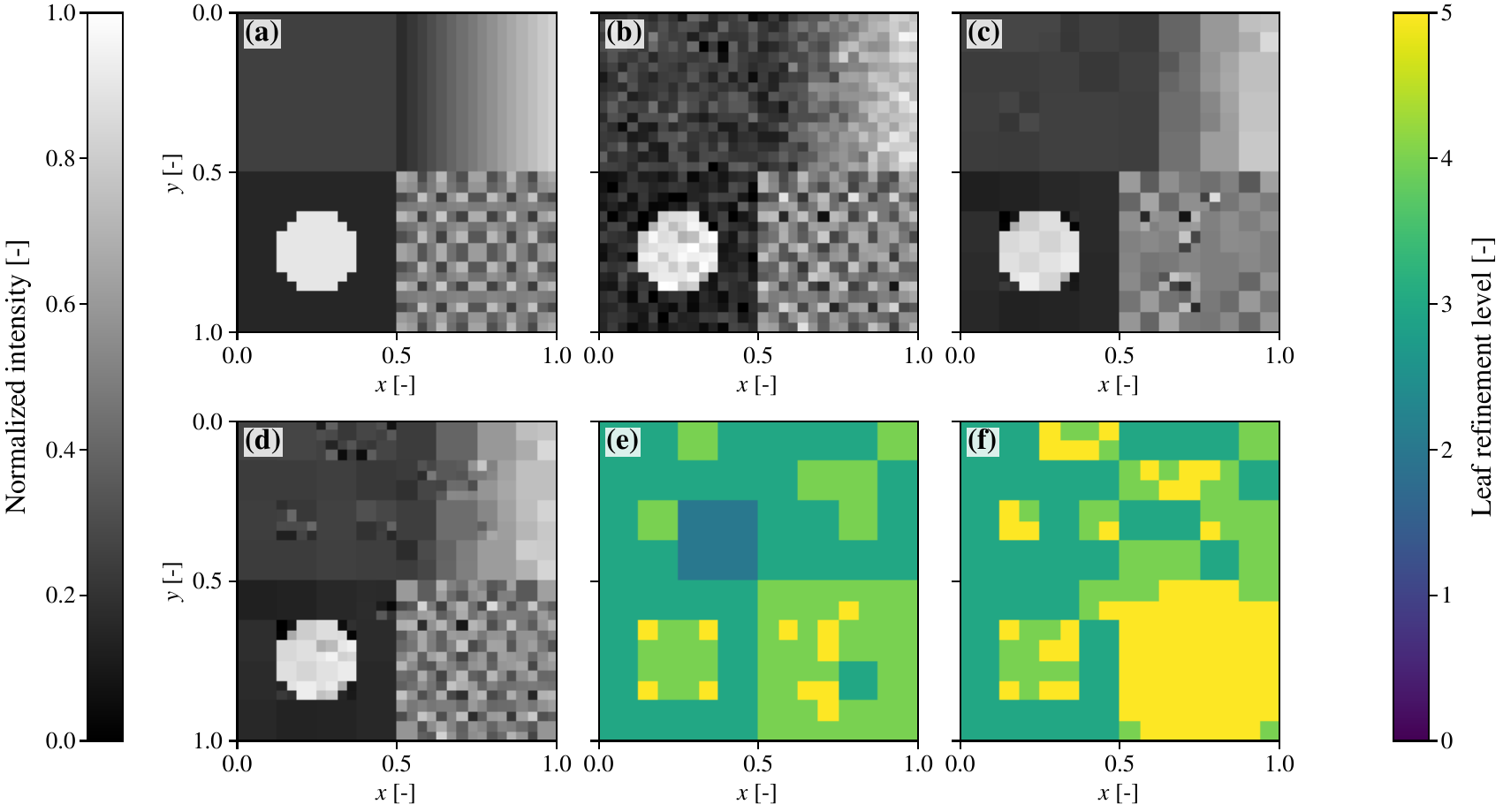}
    \caption{
    Representative spatial outcome for the four-region holdout case:
    (a) clean image; (b) noisy input; (c) deterministic-AMDI reconstruction;
    (d) Learned-AMDI reconstruction; (e) final deterministic-AMDI refinement
    levels; and (f) final learned refinement levels.  This figure illustrates
    the terminal allocations; sequential contribution is assessed separately
    by the decision-1-only control.
    }
    \label{fig:four_region_learned}
\end{figure}

\begin{figure*}[t]
    \centering
    \includegraphics[width=0.98\textwidth]
    {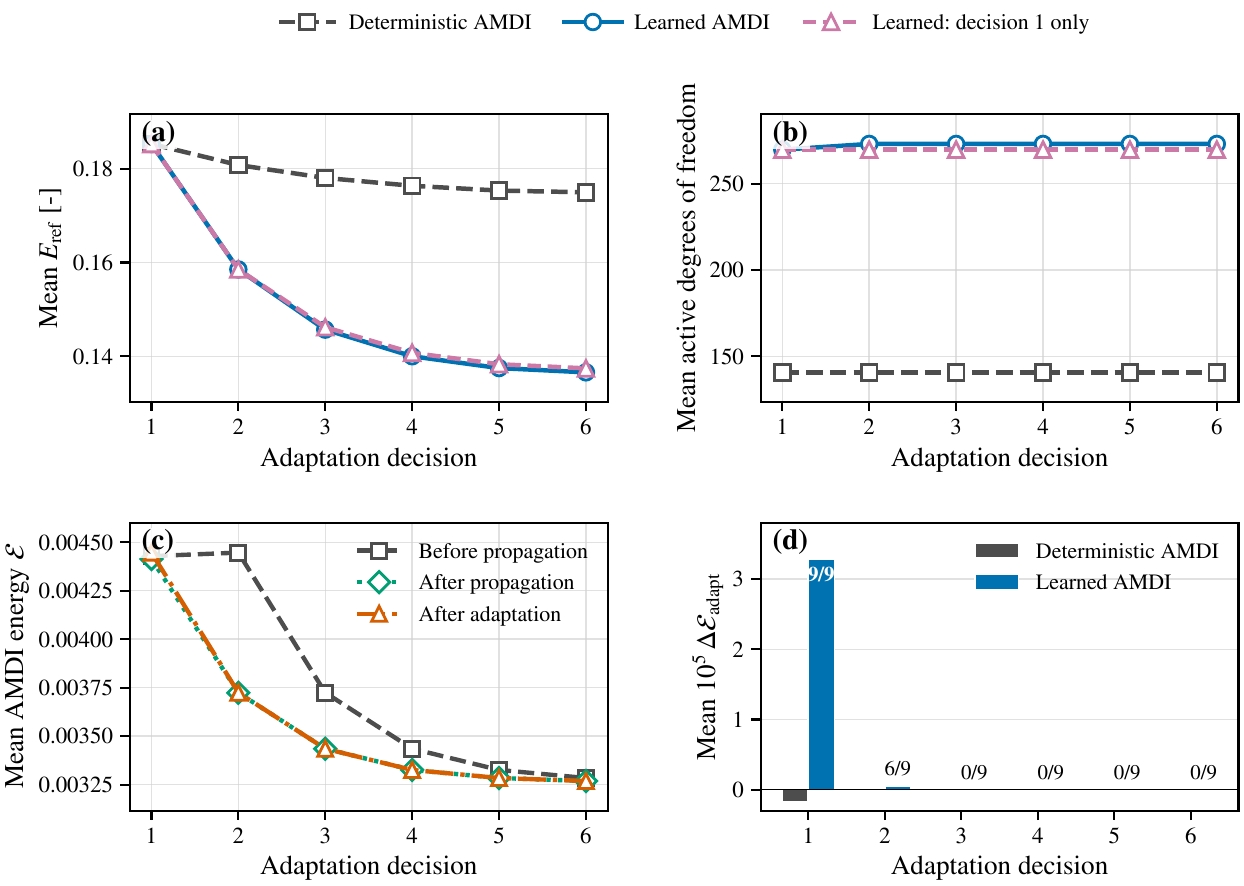}
    \caption{
    Nine-case action and energy analysis. (a) Mean reference error and
    (b) mean active degrees of freedom for deterministic AMDI, unrestricted
    Learned AMDI, and the decision-1-only temporal control. (c) Mean learned
    trajectory energy immediately before propagation, after propagation, and
    after adaptation. (d) Mean energy change caused by adaptation; labels give
    the number of cases with a positive change. The near coincidence of the two
    learned trajectories after decision~1 shows that later decisions make only a
    small contribution on this benchmark, while panels (c)--(d) show that the
    fixed-tree safeguard does not extend across the learned tree update.
    }
    \label{fig:section5_action_energy}
    \label{fig:four_region_trajectory}
\end{figure*}

%%%%%%%%%%%%%%%%%%%%%%%%%%%%%%%%%%%%%%%%%%%%%%%%%%%%%%%%%%%%%%%%%%%%%%%%%%
\subsection{Accuracy--representation tradeoff}
\label{sec:complexity}
%%%%%%%%%%%%%%%%%%%%%%%%%%%%%%%%%%%%%%%%%%%%%%%%%%%%%%%%%%%%%%%%%%%%%%%%%%

\begin{figure*}[t]
    \centering
    \includegraphics[width=0.98\textwidth]
    {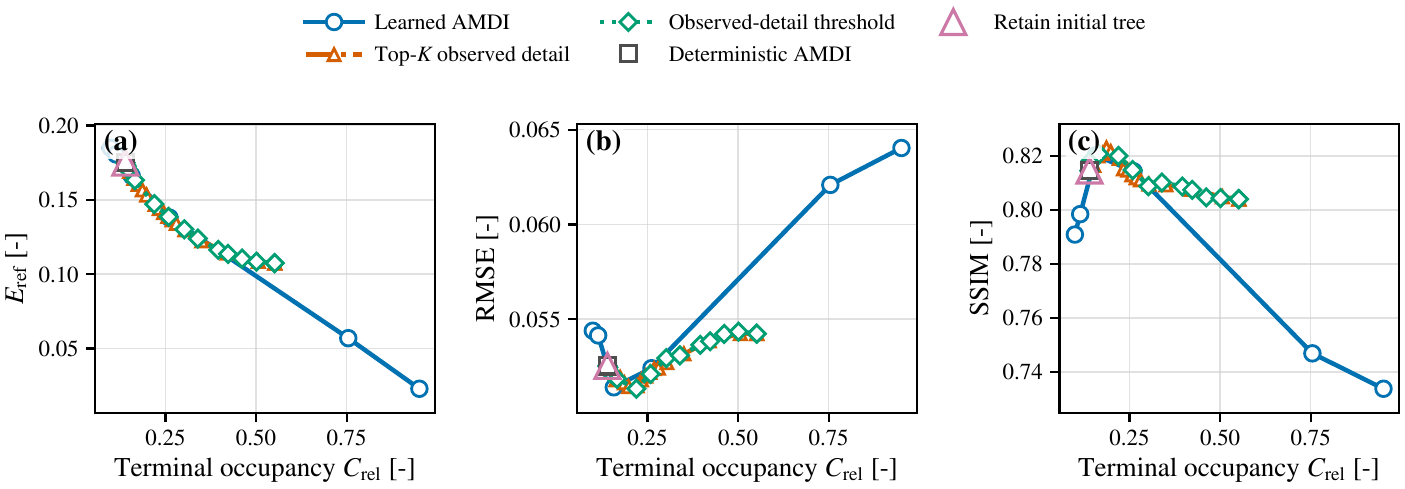}
    \caption{
    Accuracy--representation frontiers on the same four test images and noise
    realizations for every method.  Shown are (a) reference error,
    (b) RMSE, and (c) SSIM versus terminal occupancy.  Blue points are learned
    policies obtained by varying $\lambda_{\mathrm{occ}}$; the
    $\lambda_{\mathrm{occ}}=0.15$ point reuses the selected principal
    checkpoint.  Orange and green curves are the prespecified top-$K$ and
    observed-detail-threshold grids, respectively.  The square and triangle
    denote deterministic AMDI and the retain-tree control.
    }
    \label{fig:amdi2_pareto}
\end{figure*}

To examine the occupancy tradeoff, independent policies are trained with
\begin{equation}
    \lambda_{\mathrm{occ}}
    \in\{0.03,0.08,0.15,0.25,0.30,0.60\},
\end{equation}
while the remaining reward coefficients and PPO parameters are fixed.  Every
PPO point and every control point in Fig.~\ref{fig:amdi2_pareto} is evaluated
on the same four images, $s\in\mathcal{I}_4$, with noise seeds $40000+s$.
The $\lambda_{\mathrm{occ}}=0.15$ point is not a newly trained duplicate: it
reuses the selected principal checkpoint.

The nonlearned frontiers use the complete prespecified validation grids of
top-$K$ budgets and observed-detail thresholds.  The grids are fixed before
their evaluation on the four Pareto images.  Increasing the number of
selected details moves both control families smoothly from the retain-tree
point toward larger representations and smaller $E_{\mathrm{ref}}$.

The learned family spans a wider occupancy range,
\begin{equation}
    C_{\mathrm{rel}}:
    0.95093\rightarrow0.75464\rightarrow0.26099
    \rightarrow0.15698\rightarrow0.11377\rightarrow0.09912,
\end{equation}
with $E_{\mathrm{ref}}$ increasing from $0.02296$ to $0.18491$.  Within this
four-image subset, the smallest learned RMSE and largest learned SSIM occur at
$\lambda_{\mathrm{occ}}=0.25$, where
\begin{equation}
    C_{\mathrm{rel}}=0.15698,
    \quad
    E_{\mathrm{ref}}=0.16616,
    \quad
    \mathrm{RMSE}=0.05141,
    \quad
    \mathrm{SSIM}=0.82019.
\end{equation}

At the selected $\lambda_{\mathrm{occ}}=0.15$ point, Learned AMDI gives
\begin{equation}
\begin{aligned}
    (C_{\mathrm{rel}},E_{\mathrm{ref}})&=(0.26099,0.13808),\\
    (\mathrm{RMSE},\mathrm{SSIM})&=(0.05241,0.81430).
\end{aligned}
\end{equation}
On the identical cases and noise, the threshold control at
$\tau_{\mathrm{detail}}=6.12571\times10^{-3}$ gives
\begin{equation}
\begin{aligned}
    (C_{\mathrm{rel}},E_{\mathrm{ref}})&=(0.25879,0.13862),\\
    (\mathrm{RMSE},\mathrm{SSIM})&=(0.05211,0.81476).
\end{aligned}
\end{equation}
The frontiers therefore confirm the holdout conclusion: at comparable occupancy,
the validation-tuned detail rule matches the selected learned policy within
small differences and is slightly better in the image metrics.  The present
evidence supports a learned route to the allocation tradeoff, but not a claim
that PPO is intrinsically superior to observed-detail ranking on these static
synthetic images.

%%%%%%%%%%%%%%%%%%%%%%%%%%%%%%%%%%%%%%%%%%%%%%%%%%%%%%%%%%%%%%%%%%%%%%%%%%
\subsection{Resolution transfer without retraining}
\label{sec:generalization}
%%%%%%%%%%%%%%%%%%%%%%%%%%%%%%%%%%%%%%%%%%%%%%%%%%%%%%%%%%%%%%%%%%%%%%%%%%

Transfer is evaluated on four images at each resolution, using image seeds
$1000$--$1003$ and noise seed $50000+s+N$.  Table~\ref{tab:transfer_actions}
reports both the relative occupancy and the absolute number of active degrees
of freedom.  The control parameters are those selected at $32\times32$ and
are not retuned.  The retain-tree calculation coincides with deterministic
AMDI throughout this experiment because deterministic AMDI accepts no tree
change in any of the $72$ transfer decisions.

\begin{table*}[t]
\caption{
Same-resolution deployment and resolution transfer on four images at each
$N$.  $\overline{N}_{\mathrm{act}}$ is the mean absolute terminal basis size,
and $N_{\mathrm{ref}}^{(4)}$ is the total number of executed refinements over
the four trajectories.  The top-$K$ and threshold parameters are the
validation-occupancy choices $K=36$ and
$\tau_{\mathrm{detail}}=6.12571\times10^{-3}$.
}
\label{tab:transfer_actions}
\centering
\small
\setlength{\tabcolsep}{4.0pt}
\begin{tabular}{rrlrrrrrr}
\hline
$N$
& $N_{\mathrm{full}}$
& Method
& $C_{\mathrm{rel}}$
& $\overline{N}_{\mathrm{act}}$
& $N_{\mathrm{ref}}^{(4)}$
& $E_{\mathrm{ref}}$
& RMSE
& SSIM
\\
\hline
32 & 1024
& Deterministic AMDI & 0.13648 & 139.75 & 0   & 0.17301 & 0.05378 & 0.80623 \\
32 & 1024
& Learned AMDI       & 0.24927 & 255.25 & 154 & 0.13814 & 0.05466 & 0.80512 \\
32 & 1024
& Top-$K$ control    & 0.24194 & 247.75 & 144 & 0.13977 & 0.05412 & 0.81135 \\
32 & 1024
& Threshold control  & 0.24487 & 250.75 & 148 & 0.13932 & 0.05398 & 0.81207 \\
\hline
64 & 4096
& Deterministic AMDI & 0.03650 & 149.50 & 0   & 0.17534 & 0.04673 & 0.72399 \\
64 & 4096
& Learned AMDI       & 0.04089 & 167.50 & 24  & 0.17104 & 0.04404 & 0.73740 \\
64 & 4096
& Top-$K$ control    & 0.06287 & 257.50 & 144 & 0.16149 & 0.04006 & 0.78733 \\
64 & 4096
& Threshold control  & 0.04126 & 169.00 & 26  & 0.17084 & 0.04380 & 0.74187 \\
\hline
128 & 16384
& Deterministic AMDI & 0.00876 & 143.50 & 0   & 0.15712 & 0.04529 & 0.67959 \\
128 & 16384
& Learned AMDI       & 0.00931 & 152.50 & 12  & 0.15382 & 0.04318 & 0.69029 \\
128 & 16384
& Top-$K$ control    & 0.01535 & 251.50 & 144 & 0.14326 & 0.03595 & 0.74693 \\
128 & 16384
& Threshold control  & 0.00949 & 155.50 & 16  & 0.15312 & 0.04268 & 0.69274 \\
\hline
\end{tabular}
\end{table*}

The $32\times32$ transfer rows use independent noise and are therefore not
duplicates of Table~\ref{tab:holdout_main}.  Learned AMDI reduces
$E_{\mathrm{ref}}$ from $0.17301$ to $0.13814$, but its RMSE and SSIM are
nearly unchanged relative to deterministic AMDI.  The frozen threshold
control reaches essentially the same occupancy and reference error while
giving the best image metrics among the three adaptive methods in these four
same-resolution cases.

At $64\times64$ and $128\times128$, the number of learned refinements falls
from $154$ at $32\times32$ to $24$ and $12$, respectively.  The absolute
terminal basis sizes are $167.5$ and $152.5$, compared with $149.5$ and
$143.5$ for deterministic AMDI.  The frozen threshold control closely tracks
these action counts and basis sizes and slightly improves all three accuracy
measures.  By contrast, the fixed $K=36$ rule executes exactly $144$
refinements over the four cases at every resolution; its occupancy increasingly
exceeds the learned occupancy on the two finer grids and produces
correspondingly better accuracy.
It is therefore not occupancy matched after transfer.  These results show
that the actor can be deployed at finer resolutions without retraining, but
they do not show an advantage over the frozen threshold rule.

\begin{figure}[t]
    \centering
    \includegraphics[width=0.98\textwidth]
    {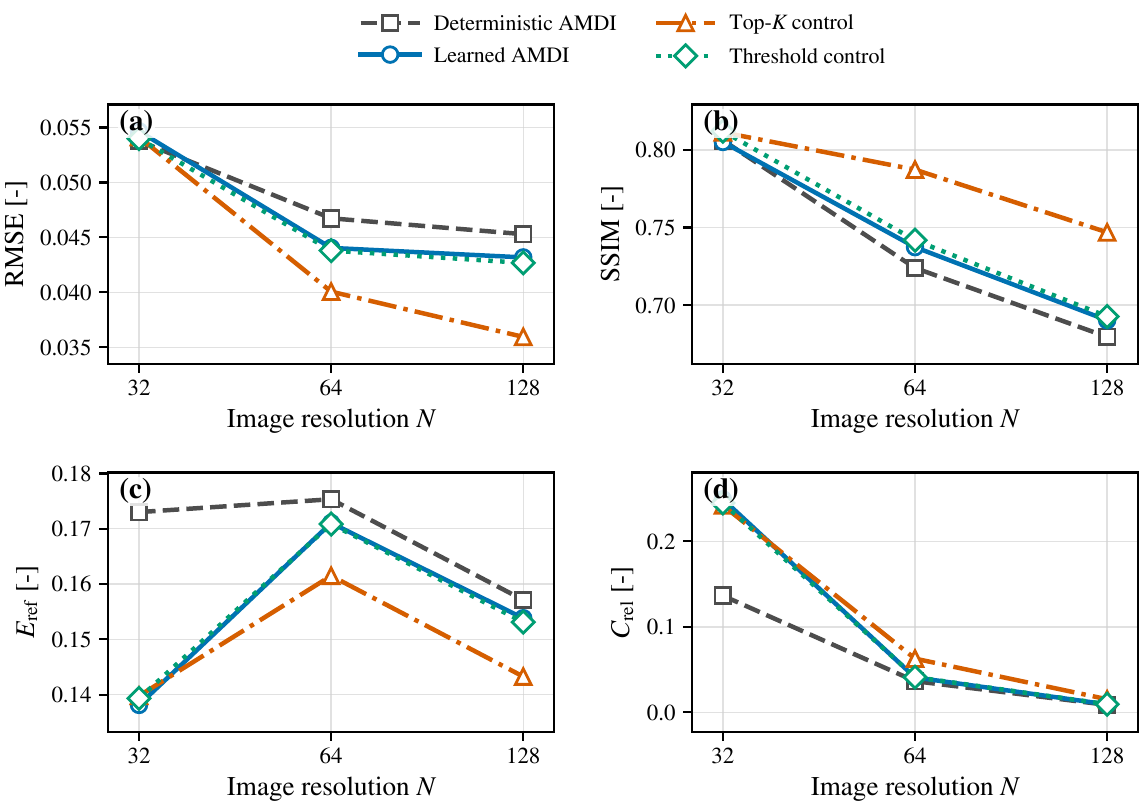}
    \caption{
    Same-resolution deployment at \(N=32\) and resolution transfer to
    \(N=64\) and \(N=128\), using four images at each resolution. Shown are
    (a) RMSE, (b) SSIM, (c) reference error, and (d) relative occupancy for
    deterministic AMDI, Learned AMDI, the frozen top-\(K\) control, and the
    frozen observed-detail-threshold control. The actor and both control
    parameters are fixed at \(N=32\) and are not retuned at the finer
    resolutions. Absolute active basis sizes and executed refinements are
    reported in Table~\ref{tab:transfer_actions}.
    }
    \label{fig:resolution_transfer}
\end{figure}

%%%%%%%%%%%%%%%%%%%%%%%%%%%%%%%%%%%%%%%%%%%%%%%%%%%%%%%%%%%%%%%%%%%%%%%%%%
\subsection{Training robustness and reward ablation}
\label{sec:robustness_ablation}
%%%%%%%%%%%%%%%%%%%%%%%%%%%%%%%%%%%%%%%%%%%%%%%%%%%%%%%%%%%%%%%%%%%%%%%%%%

Following established recommendations for reporting variation in deep
reinforcement-learning experiments
\cite{HendersonEtAl2018,AgarwalEtAl2021}, sensitivity to policy initialization
is assessed using training seeds
$20260811$, $20260817$, and $20260823$.  Each retained checkpoint is evaluated
on the same four image seeds and noise realizations.  The resulting means are
\begin{equation}
    \mathrm{RMSE}=0.05218\pm0.00183,
    \qquad
    \mathrm{SSIM}=0.81523\pm0.00907,
\end{equation}
and
\begin{equation}
    C_{\mathrm{rel}}=0.23608\pm0.05595,
    \qquad
    E_{\mathrm{ref}}=0.14886\pm0.01476,
    \qquad
    \overline{\sigma}_{\mathrm{sw}}=0.01506\pm0.00956.
\end{equation}
Here ``$\pm$'' denotes the sample standard deviation of the three
policy-seed means, not the dispersion of the twelve individual trajectories.
Each policy-seed mean is first formed from its four common evaluation cases.
Across the three training seeds, the image metrics vary less than
representation occupancy and reference error, indicating that the learned
allocation remains sensitive to policy initialization.

The reward ablation retrains after removing either the occupancy or switching
term.  Its results are summarized in Table~\ref{tab:reward_ablation}.

\begin{table}[t]
\caption{
Reward ablation on the common four-image subset.  The Full row reuses the
selected checkpoint but is evaluated with the ablation experiment's noise
realizations.
}
\label{tab:reward_ablation}
\centering
\small
\begin{tabular}{lrrrr}
\hline
Reward
& $E_{\mathrm{ref}}$
& $C_{\mathrm{rel}}$
& $\overline{\sigma}_{\mathrm{sw}}$
& RMSE
\\
\hline
Full              & 0.13769 & 0.27710 & 0.02209 & 0.05414 \\
No occupancy term & 0.02080 & 0.99561 & 0.14185 & 0.06483 \\
No switching term & 0.15126 & 0.22363 & 0.01318 & 0.05260 \\
\hline
\end{tabular}
\end{table}

Removing the occupancy penalty drives the representation to an almost uniform
state, $C_{\mathrm{rel}}=0.99561$, demonstrating that this reward component is
responsible for compression.  Removing the switching term gives a smaller
measured switching value and a smaller representation, but increases
$E_{\mathrm{ref}}$ from $0.13769$ to $0.15126$.  The switching term is
therefore part of a coupled trajectory objective; it does not control the
observed switching rate monotonically in isolation.

The Full checkpoint in Table~\ref{tab:reward_ablation}, the seed-$20260811$
checkpoint in the robustness calculation, and the principal checkpoint are
the same stored policy.  Their numerical values differ because they are
evaluated with different noise: the holdout comparison uses noise seeds
$30000+s$ (and seed $137$ for the four-region case), the ablation comparison
uses $60000+s$, and the robustness comparison uses $70000+s$.  This difference
is evaluation-noise variation, not a change in the trained policy.

\begin{figure}[t]
    \centering
    \includegraphics[width=0.98\textwidth]
    {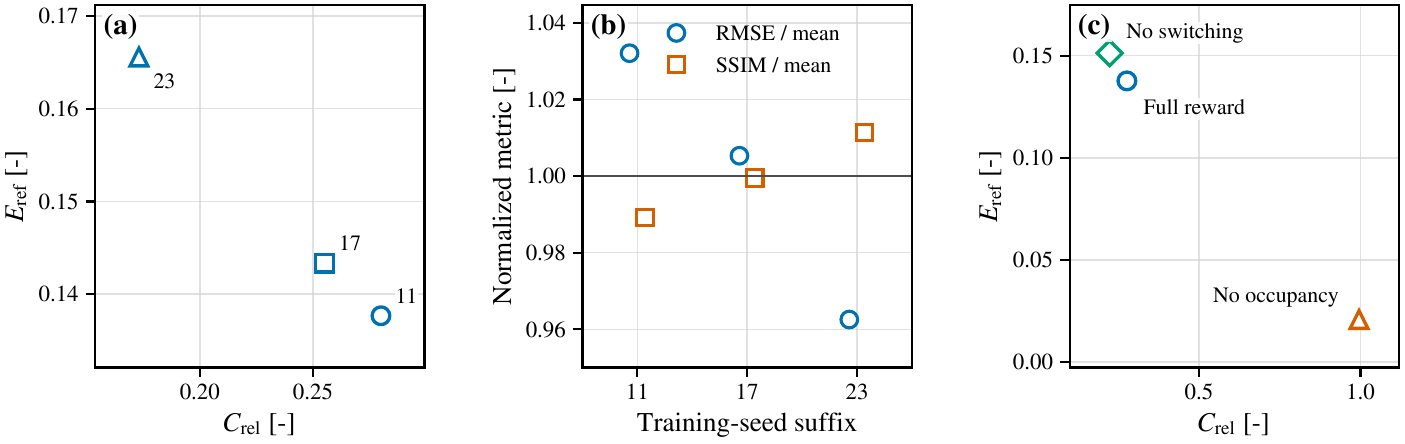}
    \caption{
    Policy-seed sensitivity and reward ablation.  (a) Error--occupancy
    operating points for three independently trained policy seeds;
    (b) RMSE and SSIM, normalized by their respective means across the three
    policy-seed means; and (c) operating points for the full reward and the two
    reward ablations.
    }
    \label{fig:robustness_ablation}
\end{figure}

%%%%%%%%%%%%%%%%%%%%%%%%%%%%%%%%%%%%%%%%%%%%%%%%%%%%%%%%%%%%%%%%%%%%%%%%%%
\subsection{Independent multiresolution localization diagnostic}
\label{sec:vampyr_crosscheck}
%%%%%%%%%%%%%%%%%%%%%%%%%%%%%%%%%%%%%%%%%%%%%%%%%%%%%%%%%%%%%%%%%%%%%%%%%%

Finally, regional localization is compared with an independent VAMPyR/MRCPP
calculation \cite{BjorgveEtAl2024VAMPyR}.  VAMPyR is used to form an order-five adaptive multiwavelet
projection of the clean analytic four-region target, with requested precision
$10^{-3}$ and maximum depth eight.  It is not run on the noisy image, the AMDI
state, or the learned trajectory.  This calculation is therefore a qualitative
localization diagnostic, not an independent validation of the learned
dynamics.

The VAMPyR effective level is \(1\) in the constant and gradient regions, has mean value \(6.54\) in the edge region, and is \(4\) in the textured region. The terminal Learned-AMDI Haar tree has regional mean levels \(3.39063\), \(3.65625\), \(3.43750\), and \(4.90625\) in the constant, gradient, edge, and textured regions, respectively. Deterministic AMDI gives \(2.87500\), \(3.25000\), \(3.31250\), and \(4.04688\). Because Haar levels and VAMPyR effective levels have different definitions, their numerical values are not compared quantitatively. The regional orderings are only partly aligned. VAMPyR assigns its largest effective level to the edge region, whereas Learned AMDI assigns its largest mean level to the textured region and does not preferentially resolve the edge over the smooth-gradient region. The comparison therefore provides a qualitative diagnostic of method-dependent localization, not evidence that VAMPyR validates the learned allocation.

\begin{figure}[t]
    \centering
    \includegraphics[width=0.72\textwidth]
    {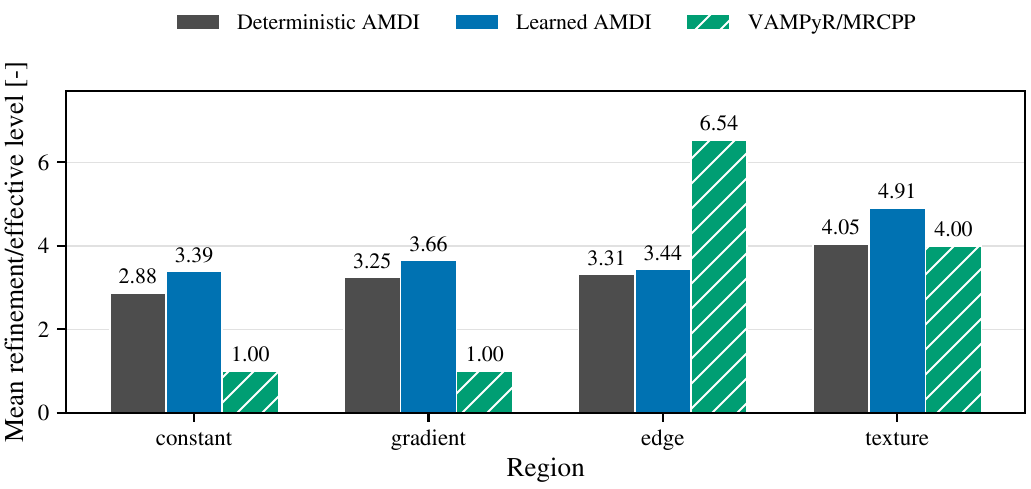}
\caption{
Regional localization for deterministic AMDI, Learned AMDI, and an
independent VAMPyR/MRCPP order-five projection of the clean analytic target,
computed with precision \(10^{-3}\) and maximum depth eight. The
method-dependent level values are used only as qualitative regional
indicators.
}
    \label{fig:vampyr_localization}
\end{figure}

\section{Conclusions}
\label{sec:conclusions}

We introduced Learned AMDI, in which a shared local policy selects the
adaptive multiresolution representation while the AMDI propagation and
hierarchy constraints remain unchanged. In the reported Haar implementation,
zero-initialized refinement candidates cannot be selected by the deterministic
one-step criterion. Accordingly, deterministic AMDI accepts no refinements in
the 54 holdout decisions and behaves as a nearly static initial-tree method.
Learned AMDI executes 393 refinements and reduces the mean terminal reference
error from \(0.17496\) to \(0.13657\), a \(21.9\%\) reduction, while
increasing terminal occupancy from \(0.13737\) to \(0.26660\). Thus, relative
to the implemented deterministic selector, the learned policy produces a
genuinely adaptive allocation and substantially improves agreement with the
uniformly resolved AMDI trajectory. The step-resolved energy data also show
that learned tree changes can introduce small energy increases, so the
fixed-tree safeguard does not imply monotonicity of the complete learned
iteration.

At nearby but not identical occupancy, the validation-tuned threshold rule
uses a slightly smaller mean representation
(\(0.26042\) versus \(0.26660\)) and gives a slightly higher reference error
(\(0.13792\) versus \(0.13657\)), while achieving slightly better RMSE and
SSIM. The learned and threshold selectors therefore occupy essentially the
same observed accuracy--occupancy tradeoff. The demonstrated role of PPO is to convert the trajectory objective into an
admissible adaptive allocation without prescribing a fixed detail threshold; the present results
do not establish superior spatial ranking at fixed occupancy.

The shared actor also transfers without retraining. Over the four transfer
images, the number of learned refinements decreases from \(154\) at
\(32\times32\) to \(24\) at \(64\times64\) and \(12\) at
\(128\times128\), adapting the absolute representation size without
resolution-specific retuning. At the two finer resolutions, Learned AMDI
improves reference error, RMSE, and SSIM relative to deterministic AMDI, while
the frozen threshold rule remains closely comparable. The decision-1-only
control shows that most of the benefit on the present static images is
produced by the initial allocation; problems with moving or emerging
localized structure provide the natural setting for testing stronger
sequential behavior. Overall, Learned AMDI provides a transparent,
hierarchy-constrained, and resolution-transferable mechanism for adaptive
multiresolution allocation.

%%%%%%%%%%%%%%%%%%%%%%%%%%%%%%%%%%%%%%%%%%%%%%%%%%%%%%%%%%%%%%%%%%%%%%%%%%
\section*{Data Availability Statement}
%%%%%%%%%%%%%%%%%%%%%%%%%%%%%%%%%%%%%%%%%%%%%%%%%%%%%%%%%%%%%%%%%%%%%%%%%%

The source code, numerical implementation, trained policy checkpoint,
experiment scripts, numerical result files, and figure-generation scripts
required to reproduce the results reported in this study are openly available
in the Learned AMDI numerical validation repository at
\url{https://github.com/Christian48596/amdi_ii_validation}.
The repository also provides the publication configuration, parameter
settings, software environment, validation tests, and instructions for
reproducing the reported calculations. All numerical data presented in this
study are generated computationally by the supplied scripts. No proprietary
or restricted datasets were used.

\bibliographystyle{siamunsrtdoi}
\bibliography{references}

\end{document}